\documentclass[11pt]{article}
\usepackage[dvips]{graphicx,color}
\usepackage{amssymb,latexsym,epsfig,epsf,xypic,amsmath}
\usepackage[mathscr]{euscript}

\makeatletter \@addtoreset{equation}{section} \makeatother

\newtheorem{Lemma}[equation]{Lemma}
\newtheorem{Theorem}[equation]{Theorem}

\newtheorem{Definition}[equation]{Definition}

\newtheorem{Conj}[equation]{Conjecture}
\newenvironment{Proof}{\noindent\emph{Proof\ }}{\hfill$\square$\\}

\newcommand\R{\mathbb R}
\newcommand\C{\mathbb C\,}
\newcommand\Z{\mathbb Z}

\newcommand\OO{\mathscr O}

\newcommand\into{\hookrightarrow}
\newcommand\res{\arrowvert_}

\newcommand\ip{{\mbox\,\_\hspace{-1.5mm}\shortmid\hspace{1.5pt}}}
\newcommand\ImO{\,\mathrm{Im}\,\Omega}
\newcommand\Jd{J\,\widetilde{\!d\theta\,}}
\newcommand\Jdh{J\,\widetilde{\!dh\,}}
\DeclareMathOperator\MCV{MCV}

\DeclareMathOperator\vol{vol}
\DeclareMathOperator\ind{ind}

\title{\textbf{Special Lagrangians, stable bundles and mean
curvature flow}}
\author{R. P. Thomas and S.-T. Yau}
\date{}

\begin{document}

\maketitle
\begin{abstract} \noindent
We make a conjecture about mean curvature flow of Lagrangian
submanifolds of Calabi-Yau manifolds, expanding on that of \cite{Th}.
We give new results about the stability condition in \cite{Th}, and
propose a Jordan-H\"older-type decomposition of (special) Lagrangians.
The main results are the uniqueness of
special Lagrangians in hamiltonian deformation classes of
Lagrangians, under mild conditions, and a proof of the conjecture
in some cases with symmetry: mean curvature flow converging to
Shapere-Vafa's examples of SLags.
\end{abstract}

\section{Introduction}

Fix a Calabi-Yau manifold $X$ with a holomorphic $(n,0)$-form
$\Omega$. In \cite{Th} a stability condition for Lagrangians in $X$
was described,
conjectured to be equivalent to the existence of a special Lagrangian
(SLag) in the hamiltonian deformation class of a fixed Lagrangian.
This was motivated by an infinite dimensional set-up
in which $U(1)$ gauge transformations act on the (infinite
dimensional) space of Lagrangians
with flat $U(1)$ connections on them. There is a natural complex
structure and symplectic form on this space and, ignoring issues
of integrability of these structures (see \cite{Th}),
the formal complexification of the $U(1)$ gauge transformations
gives hamiltonian deformations of the Lagrangian, with moment map
the $n$-form $\ImO\res L$. The stability condition was
motivated by an example of Joyce and the `angle criterion', in terms
of splittings of the Lagrangian into Seidel's graded Lagrangian
connect sums (as defined in Section \ref{index} below) and family
versions thereof, with a
certain phase inequality. This led to a conjecture, a sort of
globalised version of the angle criterion \cite{L}, \cite{N},
that the hamiltonian
deformation class of a Lagrangian should contain a SLag if and only
if the Lagrangian is stable; this SLag representative should then be
unique. Here
we expand on the conjecture and relate it to mean curvature flow. It
was verified for the simplest case $T^2$ in \cite{Th}; here we prove
it in a series of $n$-dimensional examples with symmetry (Theorem
\ref{faq}), and prove uniqueness of smooth SLags in hamiltonian deformation
classes whose Floer cohomology \cite{FO3} is defined (Theorem \ref{!}).

We write $\approx$ for ``in the same hamiltonian
deformation class as'', and use $\,\widetilde{\ }\,$ for the
isomorphism $T^*L\to TL$ induced by the metric on a Riemannian
manifold $L$. Restricting the Ricci-flat metric on a Calabi-Yau
manifold $(X,\Omega)$ to a Lagrangian submanifold $L$ we get an induced
volume form $\vol$ on $L$, and by a short calculation
\begin{equation} \label{vol}
\Omega\res L=e^{i\theta}\vol
\end{equation}
defines an $S^1$-valued function $\theta$ on $L$, the
\emph{phase function of} $L$. A \emph{grading} of $L$ is a lift of
$\theta$ to a real valued function. By Lagrangian we will always mean
graded Lagrangian (thus the Maslov class of the Lagrangian, which
is the class of $d\theta$ in $H^1(L;2\pi\mathbb{Z})$, is assumed to vanish,
and we have chosen a lift of $\theta$). $L$ is special Lagrangian
(SLag) if $\theta$ is a constant; equivalently, replacing $\Omega$
by $e^{-i\theta}\Omega$, $\ImO\res L\equiv0$.
An average, cohomological, measure of the phase of a homology class
$[L]$ is given by taking the phase of the complex number $\int_L\Omega$;
for $L$ graded with the variation of $\theta$ less than $2\pi$, this lifts
naturally to give a real number
$\phi(L)$, which is the phase of any SLag in the same homology class.

We should point out that as in \cite{Th}, we do not fully understand the
role of holomorphic discs in the theory. These are of course crucial in
the definition and hamiltonian deformation invariance of Floer
cohomology; until this is fully set up \cite{FO3} and all of its
expected properties (such as the spectral sequences of \cite{Oh2} and
\cite{P}) are proved and extended to the Calabi-Yau case,
some of the arguments below are necessarily conjectural; it will be
clear which ones. We also
deal exclusively with smooth (S)Lags; how to modify our constructions
to include singularities is an important question. Using only
(family) Lagrangian connect sums as the degenerations necessary
to describe stability of Lagrangians is also probably too restrictive,
studying other singularities and splittings may also be necessary; the
conjecture in this paper is probably just the first step in understanding
SLags in hamiltonian isotopy classes.

\noindent \textbf{Acknowledgements.} The symplectic ideas and
suggestions of Paul Seidel have been absolutely invaluable throughout
this work. We have also benefitted from comments from Kenji Fukaya,
Edward Goldstein, Spiro Karigiannis, Conan Leung, Jun Li, Elizabeth
Mann, Yong-Geun Oh and Xiao Wei Wang, and would like to thank Mike
Gage for the reference [An]. The first author is supported by a
Royal Society university research fellowship and by Imperial College,
London; the second author is supported by DOE grant DE-FG02-88ER35065
and NSF grant DMS-9803347.

\section{Mean curvature flow}

We first give a well-known geometric calculation which
we learnt from unpublished lectures of Rick Schoen on his work with
Jon Wolfson, but which dates back at least as far as \cite {HaL},
\cite{Oh1} and others.

\begin{Lemma} \label{mcf}
In the above notation, the mean curvature vector of the Lagrangian
$L\subset X$ is $MCV=\Jd$.
\end{Lemma}

\begin{Proof}
We want to show that for any vector $X$ tangent to $L$,
$X\theta=-\langle\MCV,JX\rangle$.

Picking an orthonormal basis of $T_pL$ and parallel transporting it
along rays in $L$ to a frame field $(e_i)$, $(e_i,Je_i)$ forms a local
basis for $TX$ around $p$. Letting $(f_j,g_j=-f_j\circ J)$ be the dual
basis of 1-forms, it is clear that at $p$,
$$
\Omega=e^{-i\theta}\bigwedge_j(f_j+ig_j),
$$
with $\theta$ the phase function of $L$. Since $\Omega$ is parallel,
$\nabla_X\Omega=0$ yields
\begin{eqnarray}
iX(\theta)\bigwedge_j(f_j+ig_j)&=&\sum_k(f_1+ig_1)\wedge\ldots\wedge
\nabla_X(f_k+ig_k)\wedge\ldots\wedge(f_n+ig_n) \nonumber \\
&=&\sum_k\left[\nabla_X(f_k+ig_k)\left({1\over2}(e_k-iJe_k)\right)\right]
\bigwedge_j(f_j+ig_j). \label{formula}
\end{eqnarray}

Taking covariant derivatives on the Calabi-Yau (i.e. not on $L$), we have
$$
-\langle\MCV,JX\rangle=-\langle\sum_i\nabla_{e_i}e_i,JX\rangle=
\sum_i\langle\nabla_{e_i}Je_i,X\rangle,
$$
since $J$ is both skew adjoint and parallel. But as $Je_i$ and
$X$ are orthogonal, this is
$$
-\sum_i\langle Je_i,\nabla_{e_i}X\rangle=\sum_i\langle e_i,J\nabla_Xe_i
\rangle,
$$
as we may choose $X$ to have zero bracket with the $e_i$s.

So comparing with (\ref{formula}) we are left with showing that
$$
\sum_k\left[\nabla_X(f_k+ig_k)\left({1\over2}(e_k-iJe_k)\right)\right]
=i\sum_i\langle e_i,J\nabla_Xe_i\rangle,
$$
i.e. that $(\nabla_X(f_k+ig_k))(e_k-iJe_k)=2i\langle e_k,J\nabla_Xe_k
\rangle$.

But $(f_k+ig_k)(e_k-iJe_k)=2$ is a constant, so the left hand side
is $-(f_k+ig_k)(\nabla_X(e_k-iJe_k))=-f_k(\nabla_X(-iJe_k))-
ig_k(\nabla_Xe_k)$; the other terms vanish as $\nabla_Xe_i$ was
chosen to be perpendicular to $L$. Using $\nabla_XJ=0$ and
recalling that $g_k=-f_k\circ J$, we obtain $2i\langle e_k,
J\nabla_Xe_k\rangle$.
\end{Proof}

Another simple but important result is how the phase $\theta$ and
volume form $\vol_L$ vary under a hamiltonian deformation
$J\,\widetilde{\!dh\,}$ of $L$. Such calculations appear in various forms
in \cite{Oh1}, \cite{Sm}, for instance; we give short geometric proofs for
completeness.

\begin{Lemma}
Under a hamiltonian deformation $\Jdh$ of a
Lagrangian $L$, we have
\begin{eqnarray}
{d\over dt}\theta\!&=&\!-\Delta_L(h), \label{h} \\ \nonumber
{d\over dt}\vol_L\!&=&\!-\langle d\theta,dh\rangle\vol_L.
\end{eqnarray}
\end{Lemma}

\begin{Proof}
Take real and imaginary parts of $e^{-i\theta}$
times the following:
\begin{eqnarray*}
i\dot\theta e^{i\theta}\vol_L\!\!\!&+&\!\!\!e^{i\theta}{d\over dt}\vol_L
={d\over dt}(e^{i\theta}\vol_L)=(\mathcal L_{\Jdh}\Omega)\res L=
d(\Jdh\ip\Omega)\res L
\\ &=&\!\! id(e^{i\theta}\,\widetilde{\!dh\,}\ip\vol_L)=-e^{i\theta}
d\theta\wedge(\,\widetilde{\!dh\,}\ip\vol_L)-ie^{i\theta}d^*dh\vol_L.
\end{eqnarray*}
Setting $\Delta_L=d^*d$ $(=-\sum_i\partial^2_{x_i}$ in geodesic
coordinates) gives the result.
\end{Proof}

We next show that, given a suitable metric on the Lie
algebra $C^\infty(L,\R)/\R$, the gradient flow of the norm square
$-|m|^2$ of the moment map $m=\ImO\res L$ of \cite{Th} is mean
curvature flow. The following standard calculation, applicable in any
K\"ahler reduction picture, shows that the gradient flow of
$-|m|^2$ is given by $JX_{m^*}$, where $J$ is the complex structure,
$m^*$ is the element of the Lie algebra $C^\infty(L)$ corresponding to
the moment map $m=\ImO\res L$ in the dual of the Lie algebra under the
metric on $C^\infty(L)$, and $X_{m^*}$ is its induced action on the
space \{Lagrangians with flat $U(1)$ connections on them\}.
$$
X(-|m|^2)=-X(m(m^*))=-\omega(X,X_{m^*})=\langle X,JX_{m^*}\rangle.
$$
By the definition of the group action in \cite{Th}, this
deformation $JX_{m^*}$ is just the hamiltonian deformation of the
Lagrangian $L$ with hamiltonian function $m^*$ on $L$.

Choosing the volume form Re\,$\Omega\res L$ on $L$ to define an $L^2$
metric on $C^\infty(L)$ gives $m^*=\tan\theta$, since
$m=\ImO\res L=\sin\theta\,$vol$\,=\tan\theta\,$Re\,$\Omega\res L$.
Similarly using the induced Riemannian volume form vol gives
$m^*=\sin\theta$, while using
$$
{\sin\theta\over\theta}\,\mathrm{vol}
$$
as volume form on $L$ yields $m^*=\theta$.
Any of these are suitable for small phase
$\theta:\,L\to\R$, and give similar flows down which the moment map
decreases. The last one, however, is precisely mean curvature flow, by
Lemma \ref{mcf}.

This and the previous lemma show that under mean curvature flow,
the phase $\theta$ satisfies a (time dependent) heat equation while
the Riemannian volume form decreases (as usual):
\begin{eqnarray}
\dot\theta\!&=&\!-\Delta\,\theta, \label{thetadot} \\
{d\over dt}\vol_L\!&=&\!-|d\theta|^2\vol_L. \label{voldot}
\end{eqnarray}
We therefore obtain a maximum principle for $\theta$, whose range must
always decrease, but it is important to note that the Laplacian $\Delta$
is time dependent as the metric on $L$ used to define it varies.

From these follow a series of identities and estimates, many of which
we use later, but none are strong enough to
give long term existence of the mean curvature flow, and with good
reason. Mean curvature flow is a complicated and much-studied
subject (understood only in codimension 1, dimension 1 \cite{Gr}, and,
in special cases, in
two dimensions \cite{Wa}), with known examples of finite
time blow-up. While we might expect it to behave better for Lagrangians
(locally functions of one variable instead of $n$), examples in
Section \ref{eg} show similar phenomena. But in our examples
there will be a way round these problems, and we will be able to make
a conjecture about the flow which may help in its study.

\section{Connect sums and Floer gradings} \label{index}

The stability definition in \cite{Th} made extensive use of
graded Lagrangian connect sums \cite{S2}; a description
of these and their relationship to Floer cohomology will be
important again here, as will knowledge of the Floer index of
Lagrangian intersections. We fix our conventions and definitions
now; in some places these differ in orientation from some of
the mirror symmetry literature and \cite{S2};
the problem seems to be deciding on whether to use the standard
symplectic form $dxdy$ on $T^2$, or the equally standard
$dpdq=-dxdy$ considering it as the cotangent bundle of its
SYZ base $S^1$ (divided by a lattice) \cite{SYZ}.

\subsection{The connect sum}

Suppose we have two Lagrangians $L_1,\,L_2$ hamiltonian isotoped
to intersect transversally in a finite number of points. We will
work at one of these points $p$. There we can pick a local
Darboux chart with coordinates $(x_i,y_i)$ and symplectic
form $\sum_idx_i\wedge dy_i$ such that $L_2=\{y_i=0\}$ is the
$x$-axes, and
\begin{equation} \label{def}
L_1=\{y_i=\tan(\alpha)x_i\}
\end{equation}
for some $\alpha\in(0,\pi)$. (It would be more usual to use
$\alpha=\pi/2$, of course, but that situation can be moved to
this one by an obvious symplectic (shear) transformation).

Using $z_i=x_i+iy_i$ coordinates to set up the obvious isomorphism
to $\C^n$ (notice that this complex structure and that inherited from
$X$ may be different), the $L_i$ are
$$
L=e^{i\alpha}[0,\infty).S^{n-1}:=\{z_j=re^{i\alpha}a_j\,:\,r\in
[0,\infty),\,\mathbf a=(a_j)_{j=1}^n\in S^{n-1}\subset\R^n\subset\C^n\},
$$
where $\alpha$ is set to zero to give $L_2$.

So given a curve $\gamma$ in $\C$, we define a Lagrangian
$$
L_\gamma=\gamma.S^{n-1}=\{z_j=\gamma a_j\,:\,\mathbf a=(a_j)_{j=1}^n
\in S^{n-1}\subset\R^n\subset\C^n\}.
$$
Then $L_2$ is represented by $\gamma_2=[0,\infty)\subset\C$,
$L_1$ by $\gamma_1=e^{i\alpha}[0,\infty)\subset\C$, and $L_1\cup L_2$
by the V-shaped union of these curves.

In this notation the \emph{Lagrangian connect sum} $L_1\#L_2$ is
represented by any smoothing $\gamma=:\gamma_1\#\gamma_2$
of $\gamma_1\cup\gamma_2$ staying inside the cone
$\{re^{i\beta}\,:\,r>0,\,\beta\in[0,\alpha]\}$
which is $\gamma_1\cup\gamma_2$ outside a compact set, and a
smooth curve cutting off the cone at the origin. (So here
$\gamma$ is \emph{not} a connect sum of the curves $\gamma_i$
in the topological sense; we only use the notation because the
resulting Lagrangians are topological connect sums.)

We want to analyse the phase of such a connect sum; initially
in the complex structure we picked on $T_pX\cong\C^n$ using the
$x_i+iy_i$ coordinates and so, up to scale, $\Omega\res p=
dz_1\ldots dz_n$. Then the phase function of the Lagrangian
$L_\gamma$ associated to a curve $\gamma$ is easily calculated to
be $\theta(\gamma')+(n-1)\theta(\gamma)+N\pi$ for any $N\in\Z$
(where $\theta(z)\in(-\pi,\pi]$ is the phase of a
complex number $z=re^{i\theta(z)}$). Orienting $\gamma_2$
such that $\gamma_2'$ is a positive real number, and choosing
$L_2$ to have phase 0, corresponds to choosing $N=0$ and so grading
$L_1$ by
\begin{equation} \label{pa}
(-\pi+\alpha)+(n-1)\alpha=n\alpha-\pi.
\end{equation}
In particular, choosing $\alpha=\pi/n$, and setting, for any $c>0$,
$$
\gamma^{\ }_c=
\{re^{i\theta}\,:\,r^n=c\sin(n\theta),\ \theta\in(0,\pi/n)\}
$$
gives a SLag $L_\gamma$ which has a grading of phase identically
zero, asymptotic
to $L_1$ and $L_2$ at infinity, and this is precisely the local
model of the example of Joyce, Harvey and Lawlor used so
extensively in \cite{Th}. 

While this is not strictly of the form $\gamma_1\#\gamma_2$ as defined
above (it is only asymptotic to the $\gamma_i$, not equal outside
a compact set), by taking $c$ small we can make it as close as we like
to such a connect sum, all in the same hamiltonian deformation class,
and the construction of Joyce is indeed a hamiltonian deformation of a
connect sum as claimed in \cite{Th}.

We plot these SLag curves $\gamma_c\subset\C$ in Figure 1
as the light lines, converging as $c\to0$ to the V-shaped
$\gamma_1\cup\gamma_2$ (with $\alpha=\pi/n$). Then the dark lines
depict connect sums $L_1\#L_2$ for $\phi(L_2)=0$ and $\phi(L_1)=
\pm\epsilon$. If $\phi(L_1)<0$, the \emph{stable case} as described
in \cite{Th}, then we can choose the connect sum
such that the phase of $L_1\#L_2$ varies monotonically between its
values on $L_1$ and $L_2$, i.e. between $-\epsilon$ and 0.

\begin{figure}[h]
\center{
\input{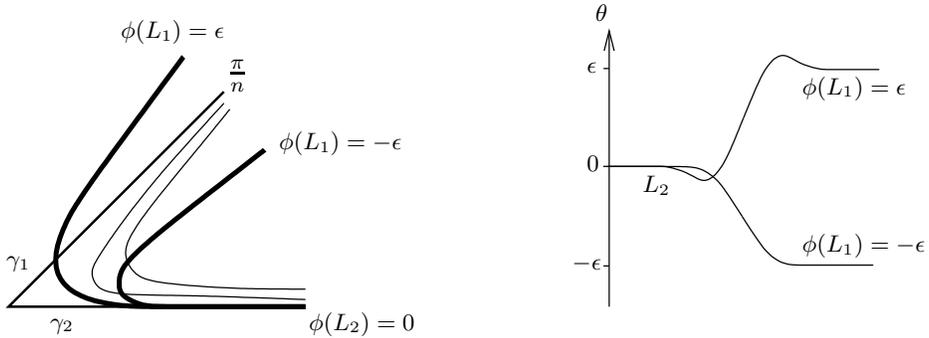}
\caption{$\gamma_1\#\gamma_2$, and the resulting phase function
$\theta_{L_1\#L_2}$, for $\phi(L_1)=\pm\epsilon$}}
\end{figure}

If, however, $\phi(L_1)>0$, the \emph{unstable case} in \cite{Th},
then the phase of $L_1\#L_2$ must initially decrease to move away
from $L_2$ before increasing to reach $L_1$ (i.e. $\gamma$ must cross
the light lines one way then the other), giving a phase
function which necessarily goes outside the range $(0,\epsilon)$
(see Figure 1).
This will be important to us later -- under mean curvature flow
we expect the phase function $\theta$ to evolve to a constant
in the stable case (under the heat equation (\ref{thetadot})) and
to a Heaviside step function (with values $0$ and $\epsilon$)
in the unstable case. This does not then contradict the maximum
principle as the unstable case has the described non-monotonic phase.

While this defines the symplectic connect sum in general by means
of our Darboux chart, the analysis of phases depended on the choice
of complex structure on $T_pX$, which may not have been the one
restricted from the Calabi-Yau $X$. In the general case we can fix
$\theta_p(L_2)=0$, without loss of generality, by rotating $\Omega$,
and pick local complex
coordinates $z_i=x_i+iy_i$ such that $\Omega_p=dz_1\ldots dz_n$ and,
at the level of tangent spaces at $p$, (the tangent space to) $L_2$
is at $y_i=0\ \forall i$. (The tangent space to) $L_1$ will be of the
form
\begin{equation} \label{l1}
L_1=\{z_i=re^{i\alpha_i}\},
\end{equation}
for some $\alpha_i$s that are \emph{no longer necessarily all the same}
in these coordinates. We are now connect summing Lagrangians of pointwise
phase $0$ and $\sum_{i=1}^n(\alpha_i)-\pi$ (compare (\ref{pa})), but
the resulting phase function will not be as simple as before -- it
is not pulled back from $\gamma$ but will vary over the $S^{n-1}$
fibres. Its average phase over the $S^{n-1}$s will have a similar
form to that in Figure 1, however, and in the case of all the
$\alpha_i$s being the same we get the earlier simpler picture.

The dependence of the hamiltonian deformation class of
$L_1\#L_2$ on the choice of scale of the neck at each intersection point
was described in (\cite{Th} Section 4) (in particular if there is
only one intersection point the class is uniquely defined). We should
also point out that the graded connect sum (when it exists)
is also independent of hamiltonian deformations of $L_1$ and $L_2$. While
the $L_i$s intersect transversely this is clear; we need only understand
what happens in crossing the codimension one wall of Lagrangians
intersecting in a double point (i.e. creating or cancelling two
intersection points). But it will be clear from the definition of grading
below that two such points must have grading differing by one, and
so the connect sum along both of them cannot be graded (\ref{grhash}).

\subsection{The grading on Floer cohomology}

The Floer cohomology group $HF^*(L_2,L_1;\C)$ \cite{FO3}
is the cohomology of a cochain complex made from a copy of $\C$
for each intersection point of two \emph{graded} Lagrangians
hamiltonian isotoped to intersect transversally.
The differential is defined by counting holomorphic
strips, with boundary in the Lagrangians, running from one intersection
point to another. It is a symplectic refinement of the
topological intersection theory of $L_1,\,L_2$, and as such is invariant
only under hamiltonian deformations of the $L_i$.
What is important to us is the grading of a particular transverse
intersection point, as defined in \cite{S2}, \cite{FO3}.

While this can be defined completely topologically, it is most easily
(and equivalently) defined via a complex structure. Again we
work at the level of tangent spaces, pick local coordinates and,
without loss of generality, take $L_2$ to have phase 0 and to be
the $x$-axes: $L_2=\{y_i=0\}$. Write $L_1$ as
$$
L_1=\{z_i=re^{i\alpha_i}\},
$$
\emph{where the $\alpha_i$s are all in} $(0,\pi)$. Then
$\sum\alpha_i=\theta_p(L_1)$ mod $\pi$ is the phase of $L_1$ up to
multiples of $\pi$, and the following integer
is the \emph{definition} of the Floer index of the point $p$:
\begin{equation} \label{maslov}
\ind_p(L_2,L_1):={1\over\pi}\left(\sum\alpha_i-\theta_p(L_1)\right).
\end{equation}
Notice therefore that $\ind_p(L_2,L_1)+\ind_p(L_1,L_2)=n$.
Applying the definition (\ref{maslov}) to the connect sums defined
in the last section (for which $\theta_p(L_1)=\sum\alpha_i-\pi$),
we recover a result of Seidel \cite{S2}:
\begin{equation} \label{grhash}
L_1\#L_2 \text{\emph{ exists as a graded connect sum if and only if }}
\ind_p(L_2,L_1)=1.
\end{equation}
(The only if
part follows from the independence of gradings and the Floer index
from the complex structure; we may therefore pick the complex
structure locally to have the form of the local model above.)
Given $L_1$ there is at most one choice of the grading on $L_2$
such that $\ind_p(L_2,L_1)=1$ at all intersection points $p$, so that
$L_1\#L_2$ can be graded.

In fact connect sums $L_1\#L_2$ whose own Floer cohomology is well
defined \cite{FO3} should correspond to Floer coclosed cochains, i.e.
elements of $HF^1(L_2,L_1)$,
mirror to the fact that extensions of sheaves $0\to E_1\to E\to E_2\to0$
correspond to elements of Ext$^1(E_2,E_1)$, as discussed in \cite{Th}.

We can also deal with the connect sums mentioned in \cite{Th}
which are \emph{relative} versions of the above construction;
$(n-r)$-dimensional connect sums carried out in a smooth family
over an $r$-dimensional base. Then the same Floer
index can be defined; there are now $r$ angles between the Lagrangians
that are zero, and $(n-r)$ whose signs can be computed to get the Floer
index. (The signs are constant over the family since the intersection
of the Lagrangians $L_1\cap L_2$ fibres over the base of the family
with fibres of constant dimension; an angle going to zero would cause a
fibre dimension to increase.)

\section{Uniqueness}

In finite dimensional symplectic quotient problems, convexity
properties of the moment map prove uniqueness of its zeros (modulo
the action of the real group) in a complexified group orbit.
Translating this into our terms is not quite possible, because there
are hamiltonian deformations of $L$ which are not given by the flow
of a \emph{fixed} hamiltonian on $L$. By this we mean $L_0,\ L_1$ are
deformations given by a constant hamiltonian $h\in C^\infty(L_0;\R)$
if the flow
\begin{equation} \label{flow}
f_t:\,L\to W, \quad {df\over dt}=\Jdh, \qquad t\in[0,1],
\end{equation}
takes $L_0=f_0(L)$ to $L_1=f_1(L)$.
All small deformations of a Lagrangian are of this form; for more
general deformations we have to use a
different proof of uniqueness of a SLag representative of a
hamiltonian deformation class (Proposition \ref{!} below), but for
these constant hamiltonian deformations
we describe the moment map proof to show how the formalism works.

\begin{Lemma} \label{?}
If two SLags $L_0,\ L_1$ are time-independent hamiltonian deformations
of each other, in the sense above, then $L_0=L_1$.
\end{Lemma}

\begin{Proof}
Without loss of generality we may take $\phi(L_0)=\phi(L_1)=0$. Then
we compute, down the flow (\ref{flow}),
$$
{d\over dt}\int_L h\ImO=\int h\,\mathcal L_{\Jdh}
\ImO=\int h\,d(\Jdh\ip\ImO)=\int\cos\theta\,dh\wedge*dh,
$$
where the last identity (equation (3.2) of \cite{Th}) is an easy
computation in local coordinates.
(We have abused notation and written $\ImO$ for $f_t^*\ImO$.)

So for $\theta$ lying in $(-\pi/2,\pi/2)$ this is always strictly
positive, but $\int_L h\ImO$ is zero at $t=0,\,1$. Thus
the two SLags must in fact coincide.

However, we must show that $\theta$ stays in this range if it starts
in it, and deal with the case when it is not so bounded. The way
to do this in fact proves the whole Lemma in one go anyway: pick a
maximum $x\in L$ of $h$. Then by (\ref{h}) $\dot\theta(x)\le0
\ \ \forall t$, but under the flow $\theta$ starts and ends at the
same value (the \emph{cohomologically} determined phase $\phi$ of the
SLags). So $\Delta h=0$ and all the second derivatives $h_{xx}\le0$
in any direction $\partial_x$ vanish. Similarly then all third
derivatives of $h$ must vanish (since we are at a maximum). Apply
the same procedure to the second derivatives $h_{xx}$ of $h$:
$\dot\theta_{xx}=-\Delta_Lh_{xx}\le0$ at $x$, since the other terms
$[\partial^2_x,\Delta_L]\,h$ involve derivatives of the metric times
third and lower order derivatives of $h$. Thus $h$'s 4th order
derivatives vanish, and so on.

To get an integral form of this, to show that
$h$ is in fact constant, it is enough to show that $h$ is constant in
a small ball around any global maximum $x$ (with $h(x)=0$, without loss
of generality). Consider geodesic balls $B_r$ of radius $r$ about $x$,
and their boundary spheres $S_r$. Fix a standard
unit-volume $(n-1)$-form $d\mu$ on the spheres $S_r$, so that the
volume form induced by the metric is $\partial_r\ip\vol=A(r)d\mu$,
where $r^{1-n}A(r)=c_1+O(r^2)$. For $r$ sufficiently small,
$d/dr\{r^{1-n}A(r)\}$ is
bounded by $c_2r$ (for some $c_2$ dependent only on the maximum of the
curvatures of $L$ in a neighbourhood of $x$ over time $t\in[0,T]$
of the flow). Therefore
$$
{d\over dr}\left\{{\int_{S_r}h\,A(r)d\mu\over r^{n-1}}\right\}=
r^{1-n}\int_{S_r}h_r\,A(r)d\mu+e,
$$
where $|e|\le c_3r^{2-n}\left|\int_{S_r}h\,A(r)d\mu\right|$.
Integrating over $r\in(0,R)$
for $R$ sufficiently small, and using the divergence theorem,
$$
f(R):=R^{1-n}\int_{S_R}h\,A(R)d\mu=
\int_0^Rr^{1-n}\left(\int_{B_r}\Delta h\vol\right)dr+E,
$$
where $|E|\le c_3\int_0^R|rf(r)|dr$.
Since $\Delta h=\dot\theta$, and $\int_0^T\dot\theta dt=0$
(where the hamiltonian deformation is over time $t\in[0,T]$,
and $\theta$ starts and ends at the same value), we see that
$$
\int_0^T|f(R)|dt\le c_3\int_0^T\int_0^R|rf(r)|dtdr
$$
for all small $R$. From this it follows that $f(r)\equiv0$.
That is, the average value of $h\le0$
over all small spheres surrounding $x$ (averaged over time
as the metric on $L$ varies) is zero, and $h$ must be identically
zero in a neighbourhhod of $x$.
\end{Proof}

However, we can do better by mirroring the algebro-geometric argument
that a non-zero map between stable bundles of the same slope is an
isomorphism, using the grading on Floer cohomology (\ref{maslov}).
This will appear to be slightly magical; the crux of the argument is
the hamiltonian isotopy invariance of Floer cohomology, provided
by precisely the holomorphic discs in the theory about which we
have had so little to say. Again we work in $n$ dimensions.

\begin{Theorem} \label{!}
Pick a connected graded Lagrangian $L$ whose obstructions
\emph{\cite{FO3}}
to the existence of its Floer cohomology vanish, and whose second
Stieffel-Whitney class $w_2$ is the restriction of a class
$\in H^2(X;\Z/2)$ on the whole manifold (for instance if $L$ is spin).

Then there can be at most one smooth special Lagrangian in the
hamiltonian deformation class of $L$.

In particular, SLag homology spheres are unique in their
hamiltonian deformation class in dimension 3 and above.
\end{Theorem}

\begin{Proof}
Since Floer cohomology is independent of hamiltonian deformations
\cite{FO3}, any two SLags $L_1,\,L_2$ in this same hamiltonian
deformation class satisfy
$$
HF^0(L_1,L_2)=H^0(L_1;\C)=\C,
$$
given that the zeroth order piece of $H^*(L)$ survives in $HF^*(L,L)$
for $L$ with Maslov class zero (\cite{FO3} Theorem E 1.7.4).
Thus there must be at least one intersection point $p$ of $L_1$ and $L_2$.

We first want to show that the (constant) phases of the $L_i$ are the same;
all we know a priori is that they differ by $r\pi$ for some $r\in\Z$.
Using a hamiltonian perturbation we may assume then that there is at least
one transverse intersection between the $L_i$s of Floer index 0, with the
phases of the $L_i$ at this point differing by $r\pi+\epsilon$. Thus writing,
locally, $L_2$ as the graph in $T^*L_1$ of $df$, $f$ has Morse index $r$
at the intersection point, so $r\ge0$. Similarly there is
a point of Floer index $n$ (i.e. a point of Floer index 0 when the roles
of $L_1$ and $L_2$ are reversed) which must correspond locally to $d$ of
a function of Morse index $n+r$; so $r\le0$ also. Therefore $r=0$ as
required.

So we have SLags $L_1,\,L_2$ of the same pointwise phase with at
least one intersection point which, if isolated, must have Floer index
(\ref{maslov}) zero. Thus $\theta_p=0$ in (\ref{maslov}), and the definition
(\ref{maslov}) of the index is therefore always positive (in fact between $0$ and
$n$), and zero only if the relative angles $\alpha_i=0\ \ \forall i$.
So the $L_i$ are tangent at $p$.

So there is no isolated transverse intersection point of Maslov class zero.
In fact, working in a small neighbourhood of the intersection, we may choose
coordinates such that $L_1$ is the graph in $T^*L_2$ of a closed one-form
$\sigma$ on $L_2$ which is also coclosed in a certain metric on $T^*L_2$
in a first order infinitesimal neighbourhood of $L_2$ (coclosedness is
the \emph{special} Lagrangian condition). In this small open set, write
$\sigma=df$, so that $d^*df=0$ and $f$ is harmonic; thus by the maximum
principle, it has no local maxima or minima. We want to show that $f$
is in fact constant. The Floer index (\ref{maslov})
of intersection points $df=0$ now reduces to the Morse index of $f$
at isolated critical points, but we also have to deal with degenerate
critical points of $f$. Assuming for a contradiction
that the critical set of $f$ is not all of $L_2$, we may perturb
$f$ inside any connected component of a small
neighbourhood of its critical set such that its value is unchanged on
the boundary, \emph{where it attains its global maximum and minimum},
and is Morse in the interior. (That we may take the extrema to be
on the boundary is the key point and a consequence of the maximum
principle.) We can then perturb $f$ further to
arrange its index 1 critical points to be lower (with respect to $f$)
than all higher index points (by general position arguments
\cite{Mi} Theorem 4.8) and then cancel any local minima with them
(\cite{Mi} Theorem 8.1). (There must be index 1 critical points if
there are any interior minima, by connectivity of our neighbourhood.)

The upshot is a hamiltonian perturbation of $L_1$, using this new
function, with no Floer index zero intersection points with $L_2$.
Thus $HF^0(L_1,L_2)=0$, a contradiction, so in fact
$f$ was locally constant and $L_1=L_2$.

The final statement follows from the fact that the obstructions of
\cite{FO3} live in $H^2(L)$, and homology spheres are spin.
\end{Proof}

As Donaldson pointed out, this proof is similar in flavour to proofs
of the Arnold conjecture. If the local situation (of all hamiltonian
deformations coming from a fixed function) held globally, the proof
would be `trivial', i.e. that of Lemma \ref{?} above. Even more
simply, if one SLag is a graph in the cotangent bundle of another,
we reduce the problem to the uniqueness of harmonic functions of
integral zero on $L$, i.e. to $H^0(L;\C)=\C$. To extend this
argument globally we need to replace de Rham cohomology $H^0(L;\C)$
by Floer cohomology $HF^0(L,L;\C)$.

\section{Analogues of some properties of sheaves}

In this section we discuss more
properties of (S)Lags that mirror those of holomorphic vector bundles
on Calabi-Yau manifolds.  As they rely heavily on Floer cohomology
arguments, many of the topics are necessarily treated informally and
unrigorously for now.

\subsection{Twisting by line bundles}

Any coherent sheaf can be twisted by a sufficiently positive
line bundle $\OO(N)$ so that it has sections; equivalently there are
homomorphisms to the bundle from any sufficiently negative line bundle.
If the sheaf has global support, this homomorphism is injective,
exhibiting $E$ as an extension
$$
0\to\OO(-N)\to E\to Q\to0.
$$
One test of our notion of subobject of Lagrangians (in terms of connect
sums), then, is that there should be appropriate connect sums mirroring
this extension.

A line bundle $\mathcal L$ defines a \emph{spherical object} \cite{ST}
of the derived category of sheaves on a Calabi-Yau manifold $X$;
that is Ext$^i(\mathcal L,\mathcal L)=H^{0,i}(X)\cong H^*(S^n;\C)$
is $\C$ in dimensions $0$
and $n$, and zero otherwise. These should be mirror to Lagrangian
homology spheres; we will consider only spheres here so that we can
use the graded Dehn twists \cite{S2} around them. Negativity
compared to some other Lagrangian may not make sense in general
(intuitively, the Lagrangian
might be mirror not to a sheaf but to an object of the
derived category with Homs in negative degrees, etc.) but instead
we can consider only those spheres $L$ with only degree zero
intersection points (\ref{maslov}) with a fixed Lagrangian $L'$.

Then it is indeed true that we can exhibit $L$ as a subobject of
$L'$: denoting by $T_L$ the (graded) symplectic Dehn twist about
$L$, simply note that
$$
L'\approx T_L^{-1}T_LL'\approx L\#[L'\#(L[\,1\,])]
$$
expresses $L'$ as a connect sum of $L$ and something else. These relations
can be shown by grading similar results in \cite{S1}. In general this will
not destabilise $L'$ due to the phase of $L$ being so negative.

\subsection{Stability of (S)Lags} \label{ssl}

It is usual in correspondences between stable objects in algebraic
geometry and solutions of the corresponding moment map PDE for one
direction of the correspondence to be reasonably straightforward to
prove, namely that objects which satisfy the PDE are stable.

While we cannot prove this for SLags, we can show, for SLags
satisfying Floer cohomology restrictions as in Theorem \ref{!}
(in particular for spheres), that they cannot be destabilised
by other SLags. (To test for stability of sheaves it is sufficient
to test only with \emph{stable} subsheaves; if the conjecture of
\cite{Th} is true then similarly we could test for stability of
Lagrangians by connect summing only SLags; this would then
be enough to prove the general stability of SLags.)

The idea is that if $\phi(L_1)>\phi(L)$, with both $L_1$ and $L$
SLags, then the Floer index of any intersection point of $L_1$
and $L$ is strictly positive (\ref{maslov}), almost by definition.
But if $L_1$ were to destabilise $L$, i.e. $L=L_1\#L_2$ for some
$L_2$, then there should be canonical morphisms $HF^0(L_1,L)\ne0$
and $HF^0(L,L_2)\ne0$, a contradiction.

The morphism from $L_1\#L_2$ to $L_2$, by which we mean an element of
\begin{equation} \label{hom}
HF^0(L_1\#L_2,L_2),
\end{equation}
can be described as follows (the element of $HF^0(L_1,L_1\#L_2)$
is similar). We use the description of the connect sum in Section
\ref{index}. Choose a Morse function $f$ on $L_2$ which has
local maxima at intersection points $p$ with $L_1$,
and in local Darboux charts as in Section \ref{index}, is
pulled up from a linear function on $\gamma_2$. Let the function
have a unique local minimum elsewhere on $L_2$, and now use
this to hamiltonian deform $L_2$ off $L_1\#L_2$. By construction
$L_1$ and $L_1\#L_2$ now intersect at the critical points of $f$
only, with Floer index the Morse index of $f$. In terms of
Figure 1, as $f$ has a maximum on $\gamma_2$ at the vertex of
$\gamma_2$, it defines a hamiltonian deformation of $\gamma_2$
\emph{downwards}, away from the connect-sum neck. As $L_2$ only
intersects $L_1$ near these connect-sum necks, we can make
our charts small enough that $L_2$ now only intersects $L_1\#L_2$
where its hamiltonian deformation intersects the old $L_2$, i.e.
at the critical points of $f$.

We now have a unique index zero point of $(L_1\#L_2)\cap L_2$ at
the unique local minimum of $f$. What we require is that
this survives in the passage to cohomology of the cochain complex
to give $HF^0(L_2,L_1\#L_2)\cong\C$. For instance, if there are no
index one points (i.e. $f$ is a Morse function with only minima and
index $\ge2$ critical points) then this will clearly be the case.
More generally there is a spectral sequence analogous to Po\'zniak's
\cite{P} with
$$
\mathrm{coker}\,\{\bigoplus_i\C_{p_i}[-n]\to H^*(L_2)\}\ 
\Longrightarrow\ HF^*(L_1\#L_2,L_2)
$$
(with a certain bigrading)
converging to $HF^*(L_1\#L_2,L_2)$. Here the notation means that a
copy of $\C$ is mapped to $H^n(L_2)$ (i.e. it is in degree $n$)
for every intersection point $p_i$ via the Morse theory for $f$ (whose
maxima are at the $p_i$). Therefore the degree zero part also survives
if, for instance, $H^1(L_2)=0$. If $L_2$ is a sphere we can proceed
more directly by applying Seidel's exact sequence \cite{S3}.

Using similar methods on Lagrangians rather than SLags, we can cut
down on the number of possible destabilising Lagrangians $L_1$
we must check to conclude that a given $L$ is stable, rather
analogously to only checking for subsheaves of vector bundles
amongst those of lower rank. There are
no morphisms (non-zero elements of $HF^0(L_1,L)$) if the phase
of $L_1$, at an intersection point $p$, is greater than that of
$L$; the Floer index at $p$ is strictly positive.
So for $L_1$ to destabilise $L$ (and so $HF^0(L_1,L)\ne0$ for
Lagrangians satisfying the same conditions as above and in (\ref{!}),
e.g. homology spheres) we must have
$$
\inf_{x\in L_1}\theta_{L_1}(x)<\sup_{x\in L}\,\theta_L(x),
$$
and in fact the corresponding phase inequality
at each point of intersection.

Thus we do not have to check all Lagrangians $L_1,\,L_2$ for the stability
of $L'$ in \cite{Th}, just those whose phase function satisfies
$$
\inf_{x\in L_1}\theta_{L_1}(x)\le\sup_{x\in L}\,\theta_L(x)
\quad\text{and}\quad
\sup_{x\in L_2}\theta_{L_2}(x)\ge\inf_{x\in L}\,\theta_L(x),
$$
where we can in fact replace the left hand sides of these inequalities
by the inf (respectively sup) over all Lagrangians in the same hamiltonian
deformation class.

Assuming the conjecture in \cite{Th}, so that we need only check
SLag destabilisers to show that stability of $L_0$, we are reduced to
checking for destabilising subobjects amongst those $L_1,\ L_2$ whose
homology classes sum to $[L_0]$ and satisfy
\begin{equation} \label{ineq}
\sup_{L\approx L_0}\left(\inf_{x\in L}\,\theta_L(x)\right)\le\phi(L_1)
\le\phi(L_2)\le\inf_{L\approx L_0}\left(\sup_{x\in L}\,\theta_L(x)\right).
\end{equation}

\subsection{A Jordan-H\"older decomposition for Lagrangians} \label{JH}

In order to understand limits of mean curvature flow it will be useful
to have the following concept;
an analogue for Lagrangians of the Jordan-H\"older filtration of
sheaves (see \cite{HuL} 1.5, for instance).

\begin{Definition}
Given two graded Lagrangians $L_1,\,L$, write $L_1\le L$ if there
exists a graded Lagrangian $L_1'$ such that $L\approx L_1\#L_1'$.
We then also write $L/L_1$ for $L'_1$, and say that $L_1$ is
a subobject of $L$.

A Jordan-H\"older filtration of $L$ is a sequence of graded
Lagrangians $L_i$ such that
$$
L_1\le L_2\le\ldots\le L_k=L,
$$
and $L_i':=L_{i+1}/L_i$ is stable. The Jordan-H\"older decomposition
of $L$ is the the singular union
\begin{equation} \label{decomp}
L_1\cup(L_2/L_1)\cup\ldots\cup(L/L_{k-1}).
\end{equation}
\end{Definition}

In sheaf theory the Jordan-H\"older filtration need not be unique,
but the decomposition is. For smooth \emph{connected} Lagrangians,
with connected $L_i$ for all $i$, however, we expect
the filtration to be unique too; the difference is essentially that
while direct sum is an operation on bundles, we are proposing that
its mirror is the (singular) union of Lagrangians, and this
cannot give a smooth Lagrangian if there is non-zero Floer cohomology
between the two Lagrangians.

If we assume the conjectures of \cite{Th}
and Section \ref{conj} below, and the properties of Floer
cohomology \cite{FO3} for all of the above Lagrangians (e.g. if they
are homology spheres), we can demonstrate the existence and uniqueness
of the Jordan-H\"older filtration for a Lagrangian $L$ whose phase
function of $L$ satisfies $\sup\theta_L-\inf\theta_L<\pi$.

Without loss of generality we may assume (by rotating $\Omega$) that 
$\theta$ lies between $\pi/2-\epsilon$ and $-\pi/2+\epsilon$, for
some $\epsilon>0$. By the inequality (\ref{ineq}) above, then, any
$L_1\#L'$ destabilising it will satisfy $\phi(L_1),\phi(L')\in
(-\pi/2+\epsilon,\pi/2-\epsilon)$.

We choose such an $L_1$ of maximal phase, and, amongst other such $L_1$s
of the same phase, minimal $\int_{L_1}\mathrm{Re}\,\Omega$ (for the
purposes of this proof we will call this quantity cohomological volume).
This still need not specify $L_1$ uniquely though.

We claim that such an $L_1$ must
be stable by construction. Any subobject of $l\le L_1$ would also be
a subobject of $L$ and so by the construction of $L_1$ must either
have smaller phase, which is not possible since it destabilises $L_1$,
or equal phase and greater or equal cohomological volume.
But $L_1=l\#l'$, where $l'=L_1/l$ has phase
$\phi(l')=\phi(L_1)\in(-\pi/2,\pi/2)$ and so positive cohomological
volume $\int_{l'}\mathrm{Re}\,\Omega$. So the complex numbers
$$
\int_{L_1}\Omega=\int_l\Omega+\int_{l'}\Omega
$$
all have positive real part, implying that the cohomological volume
of $l$ is strictly less that that of $L_1$, a contradiction.

We then apply the same
procedure to $L'$, producing an $L_2\into L'$, and so on.
By construction $\phi(L'/L_2)\le\phi(L')\le\phi(L)<\pi/2-\epsilon$,
and there is a canonical morphism (\ref{hom}) in
$HF^0(L,L'/L_2)\ne0$, making $\phi(L'/L_2)\ge\inf_L\theta_L>-
\pi/2+\epsilon$ by (\ref{ineq}).

Thus, inductively, we get the same inequalities at each stage, and
the cohomological volume of $L'$ decreases strictly with each
decomposition $L\approx L_1\#\ldots\#L_n\#L'$. The cohomological
volume of any $l$ with phase $\phi(l)\in(-\pi/2+\epsilon,\pi/2-\epsilon)$
is greater than (or equal to in the SLag case) $\cos(\pi/2-\epsilon)
\int_l\vol$, by (\ref{vol}), where $\vol$ is its
Riemannian volume form. This is bounded below above zero,
so the process can have at most a finite number of steps.

This gives us the Jordan-H\"older filtration; next we consider
uniqueness when the $L_i$s are connected (assuming the conjectures
of \cite{Th} and Section \ref{conj} and some Floer cohomology).
Suppose that $L_1'\le L_2'\le\ldots\le L$ is another such connected
decomposition where we take the $L_i'$s to be SLag assuming the
conjecture of \cite{Th}. If $HF^0(L_1',L_1)\ne0$ then
by the proof of Theorem \ref{!} (which applies as $L_1'$ and $L_1$
have the same phase), $L_1$ and $L_1'$ are equal, and we pass to $L_2$.

If, however, $HF^0(L_1',L_1)=0$, then we claim that
$HF^0(L_1',L/L_1)\ne0$. Again this should follow from standard
facts about Floer cohomology, in particular a long exact sequence
$HF^*(L_1',L_1)\to HF^*(L_1',L)\to HF^*(L_1',L/L_1)\to HF^{*+1}
(L_1',L_1)$. For $L/L_1$ a sphere this is Seidel's exact sequence
(\cite{S3} Theorem 3.3), and in general one can establish it at
the level of chains by good choices of hamiltonian perturbations
as in Section \ref{ssl}; as usual the problem is in controlling
the differential, i.e. holomorphic discs.

Assuming this we may pass to $L_2$; inductively we eventually obtain that
$L_1'$ is isomorphic to one of the graded pieces $L_{i+1}/L_i$ of
the original filtration, and is a subobject of $L_{i+1}$ but not of
$L_i$. But this gives us a contradiction (in contrast to the sheaf
analogue), since we have that both $L_{i+1}\approx L_i\#L_1'$ and
$L_1'$ is a subobject of $L_{i+1}$. The first condition ensures that
there are representatives of the hamiltonian deformation classes
such that $L_{i+1}$ and $L_1'$ have no index $n$ intersection points
by the construction of (\ref{hom}), so that $HF^n(L_{i+1},L_1')=0$.
But this is $HF^0(L_1',L_{i+1})^*$, which cannot vanish by the second
condition. (It is here we use the connectivity condition, i.e. that
the connect sum $L_{i+1}=L_1'\#(L_{i+1}/L_1')$ is not a trivial
disjoint union. Without the connectivity condition the usual
proof (e.g. \cite{HuL} 1.5) that the Jordan-H\"older
\emph{decomposition} (rather than filtration)
of sheaves is unique applies to Lagrangians, now that we have proved
or assumed all (the mirror analogues of) the algebraic facts used for
sheaves in terms of Floer cohomology instead of Exts. \\

So in the simplest case of instability, such as the example of Joyce
considered in \cite{Th}, where $L=L_1\#L_2$ is the only relevant
decomposition of $L$ with $\phi(L_1)\ge\phi(L_2)$, the Jordan-H\"older
decomposition (\ref{decomp}) would be simply $L_1\cup L_2$
(where the $L_i$ are SLag
representatives of their classes). This, like all such
decompositions, is in the closure of the hamiltonian deformation
orbit of $L$ while not being in the orbit itself.

This should have relevance to the Schoen-Wolfson programme \cite{SW} to
find canonical representatives (in a fixed hamiltonian
deformation class) of Lagrangian homology classes using volume
minimisers and so SLags (they do not use a flow, but regularity results
to study minimising currents). Our conjecture (as in \cite{Th} and later
in Section \ref{conj}) should either provide a unique SLag in a
hamiltonian deformation class, or a number of SLags in a Jordan-H\"older
decomposition.

For instance in the example above of $L_1\#L_2$ in 2 dimensions
we would produce
SLags in the classes of $L_1$ and $L_2$, but we could also form
$L_2\#L_1$; this could then be stable (it is no longer destabilised by
either of the $L_i$; if the phases of the $L_i$ are sufficiently close
one can show that in fact nothing else destabilises it either)
and we should recover a SLag in this class (and so in the same
homology class in two dimensions).

Since in two dimensions SLags are just holomorphic curves with respect
to a different complex structure, this places heavy restrictions
on stability. Take the $L_i$ above to be spheres in $K3$ surfaces.
Then any holomorphic sphere is unique in its homology class (it has
negative self intersection $-2$, so does not lie in a pencil). Any
other homologous hamiltonian deformation class
must therefore be unstable. Good examples are provided by taking a
stable (SLag/holomorphic) sphere, and applying the square of a Dehn
twist $T^2_{L_1}$ to it; this preserves homology classes but can change
hamiltonian deformation classes. If it does it should produce an
unstable Lagrangian with copies of $L_1$ in its Jordan-H\"older
decomposition; this happens in all simple cases.
$L_1\#L_2$ is taken to $L_2\#L_1$,
for instance; only one of these can be stable, the other having
a Jordan-H\"older decomposition $L_1\cup L_2$ in the simplest case.

More generally, instead of studying the action on individual (S)Lags
of symplectomorphisms like $T^2_{L_1}$ above, we could try to study
them all at once by studying the Lagrangian graph of the
symplectomorphism in $X\times X$, and its mean curvature flow. This
looks for minimal energy representatives of the hamiltonian isotopy
class of a symplectomorphism, and breaks graphs up into correspondences
representing singular maps (birational maps in the hyperk\"ahler case)
with singularities concentrated in loci whose stability is affected by
the symplectomorphism. For a Dehn twist $T_L$, for instance, we would
expect to get the graph $\Delta$ of the identity, union $L\times L$.
This also shows what the analogue of a Dehn twist $T_L$ should be when
$L$ is not a sphere but a rational homology sphere (so that it is still
spherical to complex coefficients, and so mirror to a twist on the
derived category of sheaves on the mirror Calabi-Yau \cite{ST}). Namely
$\Delta\cup L\times L$ is a Lagrangian correspondence in $X\times X$
which should give an automorphism of the derived Fukaya category of $X$
(by the usual Fourier-Mukai-type construction) \emph{not} induced by a
symplectomorphism of $X$.

\section{An example: families of affine quadrics} \label{eg}

Here we consider an example suggested to us by both Paul Seidel and
Cumrun Vafa, used in \cite{SV} and \cite{KS}. Consider the
affine algebraic variety $X^n$ given by
$$
\sum_{i=1}^nx_i^2=p(t)
$$
in $\C^n\times\C$, where $p$ is some polynomial in $t\in\C$ with only
simple zeros. Denote by $\pi:\,X^n\to\C$ the projection to the $t$
coordinate. Here we use the K\"ahler structure restricted from
$\C^{n+1}$, and the nowhere-zero holomorphic volume form given by
taking the Poincar\'e residue (\cite{GH} p 147) of the standard form
$dx_{1\ldots n}dt:=dx_1\ldots dx_{n\,}dt$ on $\C^{n+1}$;
this can be written as
\begin{equation} \label{form}
(-1)^{n+i+1}{dx_{1\ldots\hat{\imath}\ldots n\,}dt\res{X^n}\over2x_i}=
{dx_1\ldots dx_n\res{X^n}\over\dot p(t)}
\end{equation}
for any $i$ (so where $x_i=0\ \forall i$ we can use the second expression).
Here $\hat{\imath}$ means that we omit the $dx_i$ term from the wedge product.
This is then \emph{not} parallel, and the metric
we have chosen is \emph{not the Ricci-flat one}. Nonetheless it is
a good explicit testing ground for the conjecture; we can still define
$\theta$ as the phase of $\Omega\res L$ and SLags as having constant phase,
of course we then use flow
by the $\Jd$ vector, rather than mean curvature flow in
this metric. While the two flows are similar and would be the same in the
Ricci-flat metric, only the former has SLags as its stationary points
(for the latter we get minimal submanifolds, which in this metric are
not quite SLag). As Edward Goldstein pointed out to us, the $\Jd$ flow
is the gradient flow of the weighted volume functional $\int_L|\Omega|
\vol$ instead of $\int_L\vol$; everything proceeds analogously
to before on weighting all vols by $|\Omega|$, as we shall see.

Each smooth fibre over $t\in\C$ is an affine quadric with a natural
Lagrangian $S^{n-1}$
`real' slice, namely the intersection of the fibre with the slice
$$
x_i\in\sqrt{p(t)}\,\R\quad\forall i.
$$
It is invariant under the obvious $O(n)$ action on $X^n$, and is the
vanishing cycle of every singular fibre (i.e. the fibres over the
roots of $p$). Therefore any path
$\gamma:\,I\to\C$ ($I\ni u$ being some interval in $\R$) from one
zero of $p$ to another lifts to give a canonical
$O(n)$-invariant Lagrangian $n$-sphere $\gamma^n$, $S^{n-1}$-fibred
over $\gamma$ except at the endpoints where it closes up. Also,
any vector $\gamma'\partial_t$ in the base $\C\ni t$ lifts canonically
to a vector
\begin{equation} \label{lift}
\gamma'\left(\partial_t+{\dot p\over2p}\sum x_i\partial_{x_i}\right)
\end{equation}
tangent to the infinitesimal Lagrangian $\gamma^n$ lying above $\gamma'$.
Here $'$ denotes $d/du$. Note that $\gamma^1$ is a closed curve double
covering $\gamma$, branched over $\gamma$'s endpoints. We will use this
curve $\gamma^1$ later to study $\gamma^n$.

The phase function $\theta$ on $\gamma^n$ is also $O(n)$-invariant and
so a function of $t\in\C$ which we may calculate at $x_1=\sqrt{p(t)},\,
x_i=0\ \forall i\ge2$. Choosing a basis of tangent vectors to $\gamma^n$
at this point,
\begin{equation} \label{basis}
\gamma'(\partial_t+{\dot p\over2p^{1/2}}\partial_{x_1}),\,
\sqrt p\,\partial_{x_2},\ldots,\sqrt p\,\partial_{x_n},
\end{equation}
wedging them together and evaluating against the $(n,0)$-form
(\ref{form}) gives
$$
\gamma'{\dot p\over2p^{1/2}}{(\sqrt p)^{n-1}\over\dot p}=
{1\over2}\gamma'p^{n/2-1}.
$$
Therefore the phase function on $\gamma^n$ is given by
\begin{equation} \label{pf}
\theta:=\theta(\gamma^n)=\theta(\gamma')+\big({n\over2}-1\big)
\theta(p(\gamma)),
\end{equation}
where $\theta(\gamma')$ is the usual angle of the path $\gamma$, and
$\theta(p)$ is the phase of the complex number $p$ evaluated at
$t=\gamma$. So
$$
d\theta=\left({d\theta(\gamma')\over du}+\big({n\over2}-1\big)
{d\theta(p)\over du}\right)du,
$$
where $du$ is the pullback to $\gamma^n$ under the projection $\pi$
of the corresponding 1-form on $\gamma(I)\subset\C$. 

Using the metric and the orthogonal basis (\ref{basis}) we see
that
$$
\,\widetilde{\!du\,}={\gamma'(\partial_t+{\dot p\over2p^{1/2}}
\partial_{x_1})\over|\gamma'(\partial_t+{\dot p\over2p^{1/2}}
\partial_{x_1})|^2}
$$
at the point $x_1=\sqrt{p(t)},\,x_i=0\ \forall i\ge2$.

Therefore, by $O(n)$-invariance and the holomorphicity of the
projection $\pi$, $\Jd$ is the canonical lift (\ref{lift}) of
$$
{d(\theta(\gamma')+\big({n\over2}-1\big)\theta(p))\over du}\,
{\pi_*\!\left[J\gamma'(\partial_t+{\dot p\over2p^{1/2}}
\partial_{x_1})\right]\over|\gamma'|^2(1+|\dot p|^2/4|p|)}=
{{1\over|\gamma'|}{d\over du}(\theta(\gamma')+\big({n\over2}-1\big)\theta(p))
\over1+|\dot p|^2/4|p|}\,
i{\gamma'\over|\gamma'|}\partial_t.
$$
Denoting by $\mathbf t=\partial_u/|\gamma'|=\gamma'\partial_t/|\gamma'|$
and $\mathbf n=i\mathbf t$ the unit tangent and normal
vectors to $\gamma$ at a point $\gamma(u)$, the above is
$$
{\mathbf t[\theta(\gamma')+\big({n\over2}-1\big)\theta(p)]\over
1+|\dot p|^2/4|p|}\,\mathbf n.
$$
By the Cauchy-Riemann equations for the holomorphic function
$\log p=\log|p|+i\theta(p)$, $\mathbf t\theta(p)=-\mathbf n\log|p|$,
so that our flow is the lift to $X^n$ of the flow of $\gamma$ with vector
\begin{equation} \label{result1}
V^n={1\over1+|\dot p|^2/4|p|}(\MCV+\,
(1-n/2)\mathbf n(\log|p|)\,\mathbf n),
\end{equation}
where $\MCV$ is the usual mean curvature vector of $\gamma$ in
the flat metric on $\C$.

So we can reduce studying our flow to studying the flow
of a curve $\gamma$ with fixed endpoints (at zeros of $p$), under
the above vector field. We would like to relate this to mean curvature
flow of $\gamma\subset\C$ in a different metric, and also to
both our flow and the mean curvature flow for the double
$\gamma^1$ of $\gamma$ in the
double cover $X^1$ of $\C$ branched over the zeros of $p$.
The advantage of this is that we now have a flow for a closed
curve instead of a boundary value problem (but since the flow
has $O(1)=\Z/2$ symmetry it is equivalent to a flow of the
original curve $\gamma$ with fixed endpoints). We need the following
lemma.

\begin{Lemma}
Let $\langle\,.\,,\,.\,\rangle$ be the standard metric on $\C$, and
$g$ a positive real-valued function on $\C$. Then with respect
to the metric $g\langle\,.\,,\,.\,\rangle$, the mean curvature vector
of a curve $\gamma\subset\C$ is, in terms of the standard mean curvature
vector $\MCV$ (and
calculating the unit normal $\mathbf n$ in the standard metric),
$$
{1\over g}\left(\MCV-{1\over2}\mathbf n(\log g)\,\mathbf n\right).
$$
\end{Lemma}

\begin{Proof}
The endomorphism-valued 1-form $\Gamma$ defined by
$$
\Gamma_XY={1\over2g}((Xg)Y+(Yg)X-\langle X,Y\rangle\,\widetilde{\!dg\,})
$$
is symmetric and so defines a torsion-free connection on $\C$. It
is easily checked to be orthogonal with respect to the metric
$g\langle\,.\,,\,.\,\rangle$, and so gives its Levi-Civita connection
$\nabla+\Gamma$ (where $\nabla$ is the usual connection on $\C$).
Then $\langle\Gamma_{\gamma'}\gamma',\mathbf n\rangle$
(where $\mathbf n$ is calculated in the original metric) is
$-{1\over2g}|\gamma'|^2\mathbf n(g)=-{1\over2}|\gamma'|^2\mathbf n(\log g)$.

Since the unit normal to $\gamma$ in the new metric is
$g^{-1/2}\mathbf n$, the new mean curvature vector is
$$
{g\langle\gamma''+\Gamma_{\gamma'}\gamma',g^{-1/2}\mathbf n\rangle\over
g|\gamma'|^2}\,g^{-1/2}\mathbf n={1\over g}\left(\MCV-{1\over2}
\mathbf n(\log g)\,\mathbf n\right),
$$
as claimed.
\end{Proof}

Using this we can get a number of geometrically interesting flows which
are equivalent to our original flow in $X^n$. Namely, using the result
(\ref{result1}), the above Lemma, and the fact
that locally (away from branch points) $X^1$ is conformally equivalent
to $\C$ with its metric scaled by $g=1+|\dot p|^2/4|p|$ (by
holomorphicity and (\ref{lift})), we can deduce the following.

Denote by $V^n$ the flow vector $\pi_*\Jd$ of the
curve $\gamma\subset\C$ under our flow in $X^n$. Denote by
$\MCV^1_g$ the flow vector of $\gamma\subset\C$ under mean curvature
flow of $\gamma^1$ in $X^1$, with $X^1$'s natural metric scaled by
a $\Z/2$-invariant function $g$ (and omit the
$g$ in the notation if $g\equiv1$). And denote by $\MCV_g$
the mean curvature vector of $\gamma$ in $\C$ with metric $g
\langle\,.\,,\,.\,\rangle$.

Letting $\mathbf n$ be the unit normal to $\gamma$ calculated in the
standard metric on $\C$, and letting $f={|p|^{n-1}\over|p|+
|\dot p|^2/4}$\,, we have the following relations between the
various flows:
\begin{equation} \label{result2}
V^n=f.\MCV^1_f=f.\MCV_{|p|^{n-2}},
\end{equation}
and
\begin{equation} \label{result4}
V^n=\MCV^1-\left[{1\over2}|p|^{2-n}\mathbf n(f)\right]\mathbf n
=V^1-{1\over2}(n-1){\mathbf n(|p|)\over|p|+|\dot p|^2/4}\,\mathbf n.
\end{equation}

The problem with the first two is that on $\C\supset\gamma$ the flow
is \emph{not} parabolic, it has degeneracies at the end points. As
$X^1$ is so closely modelled on $X^n$ (and is in fact canonically
embedded in it), however, we might expect better on $X^1$. This is
more or less true; the result is that writing (\ref{result4})
in terms of the unit normal $\mathbf n^1$ \emph{on} $X^1$, we get

\begin{Theorem} \label{result!}
$$
V^n=\MCV^1-{1\over2}(n-1)\mathbf n^1(\log|p|)\,\mathbf n^1+
{1\over2}\mathbf n^1(\log(|p|+|\dot p|^2/4))\,\mathbf n^1.
$$
\end{Theorem}

The last term is bounded (as near a zero of $p$, $\dot p\ne0$
by nondegeneracy of $p$'s zeros) and so unimportant, we shall
see, and the flow
resulting from the first term is well understood. The second term
is more curious; it is of the order of $\mathbf t^1\theta(p)\approx
\mathbf t^1\theta((\gamma^1)')/2$ (where $'=\mathbf t^1$ denotes
differentiation with respect to arclength on $X^1$) whenever we are
close to a point where $\gamma$ emanates from a zero of $p$
(so that $p\approx Ct$ and $\theta(p)\approx\theta(\gamma)
\approx\theta(\gamma')\approx\theta((\gamma^1)')/2$).
(The last approximation is of course
not valid if $\gamma$ simply passes close to a zero of $p$; then
the equation blows up quickly as a glance at (\ref{result1}) shows,
$\gamma$ flowing to
this zero and breaking across it as discussed below; in the stable
case we will be able to rule out this behaviour and need only consider
$\gamma$ ending at the zero.)

But this is half the curvature of $\gamma^1$, so we get an
approximation to the first term again, and something like mean
curvature flow for $\gamma^1\subset X^1$. In fact in a small
neighbourhood of (the double cover of) a zero of $p$, in coordinates
$(x,y)$ in which $\gamma^1$ is a graph $y(x)$, the evolution
PDE is of the general shape
$$
y_t=y_{xx}+{y_x\over x},
$$
where by the $\mathbb Z/2$-symmetry $y_x\res{x=0}=0$, so the second
term is approximately $y_{xx}$.
So for some analysis we use this flow for $\gamma^1$, while for the rest
we pass back to $n$-dimensions, and work with the phase function
$\theta$ instead, giving a more standard (but $n$-dimensional)
parabolic equation.

We first assert how the flow behaves, before proving it in the
stable case in the next section (Theorem \ref{faq}).
Note that any deformation of
$\gamma$ is a hamiltonian deformation of $\gamma^n$ (and SLag
$\gamma^n$s have no moduli) since the $\gamma^n$s are spheres. We
picture what happens in Figures 2 and 4 in the 2 and 3
dimensional cases respectively. The dots represents zeros
of $p$ in both cases, and the epsilons and zeros are phases.

\begin{figure}[h]
\center{
\input{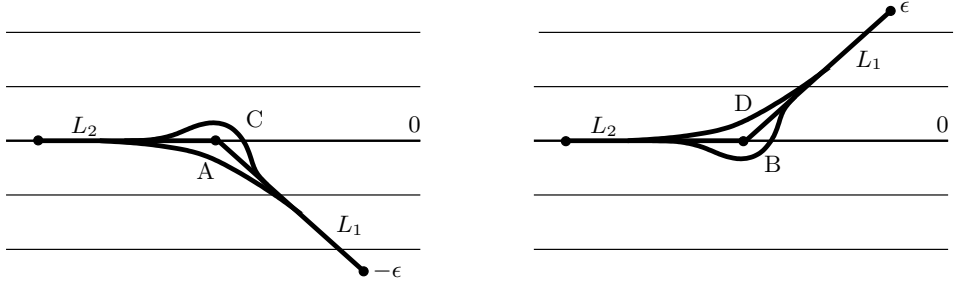}
\caption{The two connect sums $L_1\#L_2$ (1,\,2) and
$L_2\#L_1$ (3,\,4) in 2-dimensions}}
\end{figure}

In two dimensions the curves $\gamma$ whose Lagrangians $\gamma^2$ have
constant phase are the straight lines, as can be seen from (\ref{pf}).
Curves such as those marked 1 and 4 in Figure 2 flow towards
a straight line (of some \emph{non-zero} angle) corresponding to a
SLag, whereas curve 2 flows up until it `hangs' on a zero of $p$ (in
finite time), where, on restarting the flow for 2 different curves,
the separate flows form a kink and in the limit converge to
destabilising SLags $L_1,\ L_2$ of different phases.
These unions of SLags of different phases are still
stationary for the volume functional (satisfying the second
order variational equations, just not the first order SLag equations),
and in fact are minimising in odd dimensions (the angle criterion
\cite{N}, \cite{L} makes minimality of the singular union locally
equivalent to the above destabilising phase condition; in even
dimensions reversing the order of the Lagrangians reverses the
inequality and the configuration is not minimal, just stationary).

Again we see how the phase or angle criterion comes to bear; curves
1 and 2 are in the same homology class, but the two different phase
signs give very different results. As noted before in Figure 1,
this is related to the necessity of the phase to vary non-monotonically
to form the unstable connect sum; in Figure 3 we plot
the phases of the two connect sums, and with dotted lines
their limits under the heat flow (\ref{newdot}) (this is the correct
modification of (\ref{thetadot}) in the non Ricci-flat case).

\begin{figure}[h]
\center{
\input{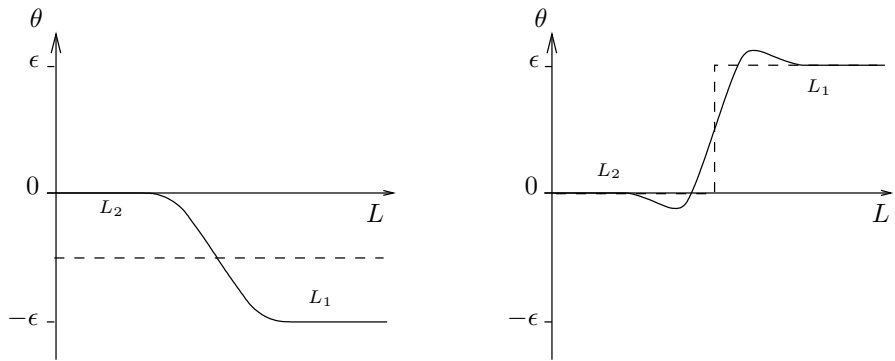}
\caption{The phase $\theta_{L=L_1\#L_2}$ of (1) $(\theta_{L_1}\equiv
-\epsilon)$ and (2) $(\theta_{L_1}\equiv\epsilon)$ respectively.}}
\end{figure}

Drawing $\gamma$ the other way round the zero of $p$ gives something in
the same homology class (the Dehn twist around the root of $p$ does not
alter the homology class of the $S^1$
fibre over $\gamma$), which is the opposite connect sum discussed in
\cite{Th} -- once the phase inequality becomes unstable for one
connect sum it becomes stable for the other.

The two connect sums are related by monodromy, as in \cite{Th}.
Take a one parameter family of polynomials $p$ which rotates
two zeros $z_1,\,z_2$ of $p$ around each other. Then under the
resulting monodromy a curve joining $z_1$ to a
third zero $z_3$ is taken from being `above' $z_2$ to
being below it, thus turning one connect sum into the other. \\

The three dimensional picture is similar. In Figure 4 we plot the
lines corresponding to SLags of phase zero, and connect sums $L_1\#L_2$
for $\phi(L_2)=0$ and $\phi(L_1)=\mp\epsilon$ (curves 1 and 2).
Again we see
the same behaviour with the phases behaving as in the graphs
in Figure 3 and the flow getting hung on a zero of $p$ in
the unstable case, splitting the Lagrangian.

Reversing the order of the connect sum in this case involves taking the
$S^2$ fibre once around the zero of $p$; this effects a Dehn twist,
reversing its orientation. Thus although curve 3 appears to give a
Lagrangian in the same homology class, it is not; the phase of $L_1$
once we have been round the root of $p$ has shifted by $\pi$ and we get
the connect sum $L_2\#(L_1[-1])$ discussed in \cite{Th}. As is also
discussed there, this can be seen to be unstable.

\begin{figure}[h]
\center{
\input{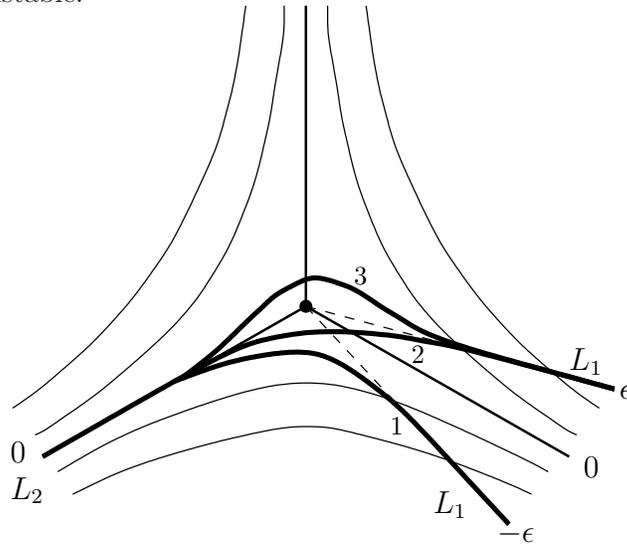}
\caption{The two connect sums $L_1\#L_2$ and
$L_2\#(L_1[-1])$ in 3-dimensions}}
\end{figure}

Things are not quite as simple as we have portrayed them if the
initial curve has very large phase variation. It is quite possible
for a curve corresponding to a \emph{stable} Lagrangian, which is
nonetheless very far from being a SLag, to pass close to a zero of
$p$ without the large negative curvature away from the zero that the
SLags exhibit. It can then flow into the zero, the Lagrangian being
split into unions of Lagrangians of which it was a connect sum,
despite their phases being such that they do not destabilise it.
This limit is in the closure of the hamiltonian deformation orbit
of the original Lagrangian, but does not contradict stability.

(A similar often-ignored subtlety occurs with stable bundles:
when moduli of semistable sheaves are created by using GIT on part of
a Quot scheme, the orbits of stable sheaves are \emph{not} closed in
Quot -- the closures contain gradeds coming from any extension of
sheaves forming the sheaf, stable or not. It is only in the
locus of points of Quot representing
semistable sheaves that the stable orbits are closed.)

The point about stable
Lagrangians is that \emph{this can be avoided} by choosing a hamiltonian
deformation of the
Lagrangian to have sufficiently small variation in phase $\theta$
(or sufficiently small volume) that it cannot be split into
destabilising Lagrangians of different phase (or higher total volume);
this we discuss now.

\section{The conjecture} \label{conj}

It is now clear what our conjecture should be. Fix a (graded) Lagrangian
submanifold $L$ of a Calabi-Yau $n$-fold $X$, and choose
the phase of $\Omega$ such that the cohomological phase $\phi(L)=0$.
Suppose first that the variation in
$L$'s phase function $\theta$ is sufficiently small in the sense that
\begin{equation} \label{close}
[\phi(L_1),\phi(L_2)]\not\subseteq(\inf_L\theta,\sup_L\theta),
\end{equation}
for all graded connect sums $L_1\#L_2\approx L$ (by this we mean
either the pointwise connect sums of Section \ref{index} or one of
the relative connect sums discussed in \cite{Th}). This condition
(\ref{close})
is preserved by the flow, by the maximum principle and equation
(\ref{thetadot}), and so prohibits $L$ splitting up as a
connect sum under the flow, as in the limit of flowing to such a
splitting (\ref{close}) would be violated.

We can also usefully consider volume instead of phase. If the
Riemannian volume of our Lagrangian $L$ is less than the
cohomological volume of any decomposition:
\begin{equation} \label{vclose}
\vol(L)\le\int_{L_1}e^{-i\phi(L_1)}\Omega+\int_{L_2}e^{-i\phi(L_2)}
\Omega,
\end{equation}
for all $L_1,\,L_2$ such that $L\approx L_1\#L_2$, then we again expect
convergence of mean curvature flow to a SLag representative for $L$.
This is also preserved under the flow, by (\ref{voldot}), and so
precludes the flow splitting $L$ into $L_1\cup L_2$, as this splitting
would necessarily have higher volume. (As the referee pointed out,
this is a global condition on $L$, whereas (\ref{close}) appears
to be a pointwise one. The global nature of (\ref{close}) is that the
points concerned are the extrema of $\theta$ over all of $L$.)

\begin{Conj}
If $L$ satisfies either of the conditions (\ref{close}) or
(\ref{vclose}) then mean curvature flow for $L$ exists for all time
and converges to a special Lagrangian in its hamiltonian deformation
class; the unique SLag conjectured in \emph{\cite{Th}}.
\end{Conj}

It is of course a consequence of this and the conjecture in \cite{Th}
that some hamiltonian deformation of $L$ satisfies (\ref{close}) if
and only if it is stable. The SLag should also be unique in its
hamiltonian deformation class as in Theorem \ref{!}.
If $L$ is stable but not close enough to
being SLag that (\ref{close}) fails, then mean curvature flow
can become singular in finite time, (locally) splitting the Lagrangian
in the reverse of a connect sum operation (i.e. with a vanishing cycle
which is an $S^{n-1}$, or an $S^{n-r}$-bundle over an
$(r-1)$-dimensional base in the relative connect sum case). We might
then conjecture that the resulting pieces are smooth so we can
begin the process again until we get a
decomposition into different phase SLags. Typically, in the simplest
case, we would get $L_1\cup L_2$ (with $\phi(L_1)<\phi(L_2)$ by
stability) which is not a hamiltonian deformation of $L$ (though it
is in the closure of such deformations).

If $L$ is unstable, we would again expect such finite time
singularities and SLag splittings. But if $L$'s phase variation,
or volume, is sufficiently small, we can hope for convergence
to the Jordan-H\"older decomposition of Section \ref{JH}. That is,
while the volume of $L$ must be larger than the cohomological volume
of its Jordan-H\"older decomposition, if it is less than any other
decomposition then it can only flow to the former. Again we expect
the flow to become singular in finite time, the limit
(locally) splitting $L$ into pieces for which we restart the flow.
This splitting of the Lagrangian is a manifestation of the well
known finite-time dumb-bell singularities in mean curvature flow.

\subsection*{Proof for our example}

Under our flow (\ref{result!}) in the Shapere-Vafa example
the proofs of the evolution equations for the phase
function $\theta$ (\ref{thetadot}) and volume (\ref{voldot})
show that the equations are modified (as the metric is not quite
Ricci-flat) to
\begin{eqnarray}
{d\over dt}\theta\!&=&\!-\Delta\,\theta+{\langle d\theta,d|\Omega|
\rangle\over|\Omega|}, \label{newdot} \\
{d\over dt}(|\Omega|\vol_L)\!&=&\!-|d\theta|^2(|\Omega|\vol_L).
\label{nvdot}
\end{eqnarray}
when $|\Omega|=|\Omega/\!\vol\!|$ is not $\equiv1$.
Therefore the maximum principle still holds for $\theta$, and
the condition (\ref{close}) is again preserved by the flow. Similarly
if we measure volume with respect to $|\Omega|\vol_L$ then this
is decreasing and (\ref{vclose}) is preserved by the flow. We can
now prove the appropriate version of our conjecture in this example.
From the proof it will also be clear that the original conjecture
could be proved in this case in the $O(n)$-invariant Ricci flat
metric if we knew it,
we would just not be able to be as explicit about the flow equations.

\begin{Theorem} \label{faq}
Suppose that $\gamma$ is a curve in $\C$, with endpoints at zeros of
$p$, and otherwise missing the zeros of $p$, such that its
pointwise phase $\theta$ (\ref{pf}) satisfies (\ref{close})
for all Lagrangians $L_i=\gamma^n_i,\ i=1,2,$ fibred over curves
$\gamma_i$ in the base, and also $S-I:=\sup_\gamma\theta-
\inf_\gamma\theta<2\pi/3$. Then the flow (\ref{result!}) exists for all
time and converges in $C^\infty$ to a smooth curve whose phase
function (\ref{pf}) is constant.
\end{Theorem}

We break the proof up into existence of the flow (best dealt with
at the level of $\gamma^n\subset X^n$), controlling the angle
variation (using $\gamma\subset\C$) to ensure no $180^o$ kinks
appear in $\gamma$, and using this to show the
flow exists for all time (for which we use the double cover
$\gamma^1\subset X^1$
and $\theta$ on $\gamma^n$). We follow \cite{An}, in parts heavily
modified to take care of the endpoints of $\gamma$. Finally
we will show that the flow converges to a SLag.

\begin{Lemma} \label{exist}
The flow (\ref{result!}) exists while the curvature of $\gamma^1$
is bounded.
\end{Lemma}

\begin{Proof}
Firstly, short term existence of the flow, given any initial curve
$\gamma\subset\C$ missing the zeros of $p$ except at its endpoints
and \emph{such that $\gamma^1$ is $H^{2+\alpha}$} for some
$\alpha>0$ (i.e. $\gamma^1$ has H\"older continuous curvature),
is in fact most easily proved at the level of the $H^{2+\alpha}$
Lagrangian $\gamma^n$; see \cite{An} for the method in 1 dimension
(which easily generalises to $n$ dimensions), and \cite{Ch} for
a similar $n$-dimensional result. This is also done
in (\cite{Sm} Proposition 1.6) using results of Hamilton \cite{H},
for instance.

While the curvature of the curve $\gamma^1\subset X^1$ is bounded,
so is the norm of the flow vector (the last term in (\ref{result!})
is always bounded, and the second term can be bounded by the
curvature at an intermediate point by Taylor's theorem). So at
any finite time $T$ the flow converges to a limit curve $\gamma_T^1$
pointwise. Parametrising the curves by their arclength on $X^1$,
their first and second derivatives as maps to $X^1$ are therefore
bounded, which
by Arzel\`a-Ascoli implies that for a subsequence of $t$ we have
convergence in $C^1$ to a $C^1$ curve with bounded (weak)
curvature. By the uniqueness of the limit, then, $\gamma_T^1\subset
X^1$ has bounded curvature.

Bounds on (the derivative of) the phase of $\gamma^1$ give
corresponding bounds on (the derivative of) the phase of
$\gamma^n$ (via (\ref{result!}) for $n$ and $n=1$). So the phase
function $\theta$ of $\gamma^n$ is also $C^0$ convergent to the
phase of $\gamma_T^n$, and satisfies the parabolic equation
(\ref{newdot}). Putting this into local coordinates and
differentiating with respect to arclength $s$, we get a
uniformly parabolic equation with bounded coefficients and a
bounded solution $\theta_s$ on $t\in[0,T]$. By (\cite{LSU}
Section III Theorem 10.1), then, $\theta_s$ is
in fact $\alpha$-H\"older continuous for some $\alpha>0$,
and $\gamma^1_T$ is $H^{2+\alpha}$. By the existence of the flow
for $H^{2+\alpha}$ initial conditions, then, the flow
exists for some time $t>T$.
\end{Proof}

\begin{Lemma} \label{cone}
$\limsup_{|s-s'|\to0}|\theta(\gamma'(s,t))-\theta(\gamma'(s',t))|<\pi$
for all time $t$ for which the flow exists, where $s$ is arclength
along $\gamma(\,.\,,t)$, and $\theta(\gamma')$ is the argument
of $\gamma'\in\C$.
\end{Lemma}

\begin{Proof}
Working outside a fixed neighbourhood of the zeros of $p$ at the
endpoints of $\gamma$, this follows from the bounds
on $\theta=\theta(\gamma')+(n/2-1)\theta(p(\gamma))$ (\ref{pf})
coming from the $\sup\theta-\inf\theta<2\pi/3$ assumption
(preserved under the flow by the maximum principle (\ref{newdot})),
as the variation of $\theta(p(\gamma))$ can be made arbitrarily small
with $|s-s'|$. Since $\gamma$ must stay at a bounded distance from
other zeros of $p$ by the condition (\ref{close}) and the
maximum principle for $\theta$ (\ref{newdot}), we are left with
proving the lemma in an arbitrarily small neighbourhood in $\C$
of the endpoints of $\gamma$.

Unfortunately, in this region, the bounds we want for
$\theta(\gamma')$ do not follow from the
bounds we have for $\theta$, and in fact come only from
comparison with known solutions. Draw the SLags of phase $S,\,I$
emanating from a zero of $p$, i.e. the curves in $\C$ solving
$\theta(\gamma')+(n/2-1)\theta(p(\gamma))=S$ or $I$.
(Since this is an ODE, there is no problem in finding solutions
and extending them to either infinity or another zero of $p$; see
\cite{SV}.) In the tangent space to the zero of $p$ this gives a
cone of angle $(S-I)/n<2\pi/3n$ which $\gamma$ lies inside and
cannot cross either at $t=0$ or any later time in the flow. So in
a sufficiently small neighbourhood of an endpoint of $\gamma$,
we may bound $\theta(\gamma)$ inside a cone of angle less than
$2\pi/3n$, and also take $\theta(p(\gamma))$ to be within any given
$\epsilon$ of $\theta(\gamma)+C$ (since the zero
of $p$ is nondegenerate; here $C$ is the phase of $\dot p$ at the
zero of $p$). Thus $\theta(p(\gamma))$ can be
bounded inside a similar cone, so that (\ref{pf}) bounds the
variation of $\theta(\gamma')$ by $2\pi/3+(n/2-1)2\pi/3n<\pi$.

The bounds on $\theta(\gamma')$ imply that the curve does not
spiral round its endpoints but moves away from them with nonzero
derivative inside the above cone until it is outside the small
neighbourhood employed above. So the remaining case to consider
is if the curve can pass arbitrarily close to one of its own
endpoints at some bounded-below arclength from its endpoints,
i.e. if the cone of SLags above starting from a zero of $p$
passes either side of that same zero at some nonzero arclength.
Then there would be a Slag fibred over a curve starting and
ending at the same root of $p$, with our Lagrangian $\gamma^n$
the connect sum of this SLag and some other $\gamma_1^n$.
But $\gamma^n$ cannot flow arbitrarily close to such a connect sum
as its phase variation would approach at least
the difference between the phase of the SLag and $\phi(\gamma^n_1)$,
contradicting (\ref{close}).
\end{Proof}

By Lemma \ref{exist} the flow exists for all time unless, as we
suppose now, the curvature of $\gamma^1$ becomes unbounded in finite
time. Then to get a contradiction we start by scaling as in \cite{An}.
Pick $s_i,\,t_i,\ i=1,2,\ldots$ such
that $t_i$ tends to the blow up time and the curvature $\kappa_i$ of
$\gamma^1(s_i,t_i)=y_i$ is maximal
over the curvatures of $\gamma^1(s,t)$ for all $s$, and all
$t\le t_i$ (here we parameterise by arclength $s$ on $X^1$, centred
at a zero $z$ of $p$, i.e. $\gamma^1\res{s=0}=z$ lies over an
endpoint of $\gamma$).

How we handle the blow up depends on
whether it happens at the branch points of $\gamma^1$ (i.e.
the endpoints of $\gamma$), by
which we mean $|s_i|=O(|\kappa_i^{-1}|)$, or in the interior
$|s_i|\gg|\kappa_i^{-1}|$, due to the different nature of
(\ref{result!}) at the branch points.
We first deal with the interior where the flow is a perturbation
of mean curvature flow and so can be handled by \cite{An}:

\begin{Lemma}
Supposing that the curvature blows up as above, then
$|s_i\kappa_i|$ is bounded.
\end{Lemma}

\begin{Proof}
Firstly, if after passing to a subsequence of $i\in\mathbb N$,
and centring $s$ about the other branch point of $\gamma^1$
if necessary, the blow up
occurs at a finite distance $|s_i|>\epsilon>0$ from either branch
point of $\gamma^1$ 
then in this interior the flow (\ref{result!}) is a finite
perturbation of mean curvature
flow satisfying the conditions of \cite{An}, so a $180^o$ kink
must appear in $\gamma^1$, contradicting Lemma \ref{cone}. So
we need only deal with the case of $s_i\to0$ (by passing to a
subsequence to concentrate around one of the two branch
points, if necessary) while $r_i:=|s_i\kappa_i|\to\infty$.

Then we rescale as in \cite{An};
\begin{equation} \label{scale}
s\mapsto \kappa_is, \quad g\mapsto \kappa_ig, \quad
t\mapsto \kappa_i^2(t-t_i),
\end{equation}
where $g$ is the metric on $X^1$. This rescaled flow for
$\gamma^1_i$ has the same form as (\ref{result!}),
$$
\dot\gamma_i^1=\MCV-(n-1)\mathbf n(\log|p|_i^{1/2})\,\mathbf
n+{1\over2}\mathbf n(\log(|p|_i+|\dot p|_i^2/4))\,\mathbf n,
$$
but with $|p|$ and $|\dot p|$ replaced by their pullbacks
$|p|_i,\ |\dot p|_i$ to
the new Riemannian surface (here it is important that $\dot p$
is still computed in the old coordinates, \emph{then} pulled
back). Therefore their gradients are scaled by $\kappa_i^{-1}$.
$\mathbf n$ denotes the unit normal to $\gamma^1$ in the
new metric on $X^1$. The curvature gets scaled by
$\kappa_i^{-1}$ and so has a maximum, over $t\le0$, of 1 at
$y_i$ (at time $t=0$). We want to show that the two perturbation
terms on the right hand side of the above flow tend to zero as
$i\to\infty$.

In the rescaled variables, work in a geodesic disc in $X^1$ of
radius $r_i/2$ ($r_i=|\kappa_is_i|\to\infty$ as $i\to
\infty$) about $y_i$. As $s_i\to0$, for
$i$ sufficiently large this is within an arbitrarily small
neighbourhood of $z$ in the original metric,
in which $\gamma'$ varies within an angle $<\pi$ cone
as in the proof of Lemma \ref{cone}, i.e. $(\gamma^1)'$
varies within an angle $<\pi/2$ cone on the double cover
$X^1$. So arclength $s$ on $\gamma^1$ and radial distance $r$
in $X^1$ are equivalent metrics on $\gamma^1$ in this disc;
$r_s:=\partial r/\partial s$ and $s_r=\partial s/\partial r$
are both bounded.

As $y_i\in\gamma_i^1$ is of arclength $s_i\ge r_i$ from $z$ (the
zero of $p$) at $s=0$, we deduce that all points of our disc
are of distance $\ge cr_i/2$ from the zero of $p$ (for some
constant $c>0$ fixed for all $i\gg1$) in the new metric.
Thus, for $i$ large enough, we have
$$
|p|_i^{1/2}\ge C\kappa_i^{-1}(cr_i/2),
$$
where $C$ is a constant just less than
the norm of the derivative of $p^{1/2}$ at the zero $z$
in the original metric on $X^1$ ($p^{1/2}$ pulls back
to a well defined function on $X^1$ with a simple zero
at $z$). We can therefore bound
$$
|\dot\gamma_i^1-\MCV|\le(n-1){\kappa_i^{-1}\sup|d(p^{1/2})|
\over C\kappa_i^{-1}(cr_i/2)}+{1\over2}\kappa_i^{-1}
\sup|d\log(|p|+|\dot p|^2/4)|,
$$
where both sups are taken over small neighbourhoods of
$z$ \emph{in the original metric on} $X^1$. As $i\to\infty$,
$\kappa_i,\,r_i\to\infty$, so the above bound tends to zero,
while the radius of the disc we are working on
$r_i/2\to\infty$. It follows that in the limit we get
mean curvature flow of a curve inside an infinite flat
disc $\R^2$; see (\cite{An} Section 9) for how to pass to the
limit to conclude that for this
blow up to occur a $180^o$ kink must appear in the curve
$\gamma^1$ (by which we mean the limsup in Lemma \ref{cone}
is $\ge\pi$). But this contradicts Lemma \ref{cone}.

(We do not repeat Angenent's argument here as we will
give a slightly harder, $n$-dimensional, version of
it around the endpoints of $\gamma$ in Lemma \ref{blowup}
below. The point is just that in the rescaling we can
get rid of the last two terms of our flow to reduce to the
results of \cite{An}.)
\end{Proof}

The remaining case
we must dismiss is that of $\kappa_i$ blowing up at points
$y_i=\gamma^1(s_i,t_i)$ with $|s_i|<A/|\kappa_i|$ for some
fixed $A$. Here we must work harder than \cite{An}.

\begin{Lemma} \label{blowup}
The curvature of $\gamma^1$ does not blow up in finite time.
\end{Lemma}

\begin{Proof}
By Lemma \ref{cone} we know that there is an $A$ such that
$|s_i|<A/|\kappa_i|$.
We rescale variables as in (\ref{scale}), and work on
a length $\kappa_i^{1/2}\to\infty$ interval (in the new metric)
on $\gamma^1$ centred ($s=0$) at the zero $z$ of $p$. This is
contained inside the ball of radius $\kappa_i^{-1/2}\to0$
about $z$ in $X^1$ in the original metric, so for $i$ sufficiently
large we can assume that $p(t)-Ct$ is arbitrarily small
\emph{in $C^2$ norm} (here $t\in\C$ is the base parameter,
\emph{not} time, $C=\dot p(z)$, and the same is true of any $C^r$
norm; $r=2$ is the case of interest for us). We start by obtaining
bounds on the polar angle of the curve and its tangent vector.
We shall confuse functions on $\C$ with their pullbacks to $X^1$
(so writing things like $p(\gamma^1)$ etc.).

Taking $i$ sufficiently large that the metric on the
radius $\kappa_i^{1/2}$ disc about $z$ in $X^1$ is sufficiently
close to being flat, define geodesic polar coordinates on $X^1$,
$r^1:=|\gamma^1|,\ \theta^1:=\theta(\gamma^1)$ (which is
$\theta(\gamma)/2$ to within a constant). Then we can
assume that $\theta^1_s$ is arbitrarily $C^1$ close to
\begin{equation} \label{sin}
{1\over r}\sin(\theta(\gamma^1_s)-\theta^1),
\end{equation}
which is the exact formula for a flat metric and polar
coordinates. This bounds $|r\theta^1_s|$. Since by
construction the curvature of $\gamma^1$ is not more than
one, i.e. $|(\theta(\gamma^1_s))_s|\le1$, we can bound
$|\theta(\gamma^1_s)|\le s$.

Note that (\ref{sin}) in flat space gives us the differential equation
$$
f_s=\kappa_i-{\sin f\over r}
$$
for $f=\theta(\gamma^1_s)-\theta^1$, with $f(0)=0$ and $|\kappa_i|\le1$.
This implies that $|f(s)|\le|s|$ (consider a point where the graph
of $f$ crosses that of $\pm s$, where $|f_s|\ge1$, for a contradiction),
so for $i$ sufficiently large that our polar coordinates are sufficiently
close to flat coordinates we can deduce a bound on $|\theta(\gamma^1_s)-
\theta^1|/s$. Thus $\theta^1/s$ is also bounded, and $\theta^1/r$ by
the uniform comparison bounds of $r$ and $s$ given by the cone argument
in Lemma \ref{cone}.

Instead of considering the equation (\ref{result!}) for $\gamma^1$,
we analyse the equation (\ref{newdot}) for $\theta$ on $\gamma^n$.
After rescaling it becomes
\begin{equation} \label{eqn}
\dot\theta=\Delta\theta+{\langle d\theta,d(|\Omega|_i)\rangle\over
|\Omega|_i},
\end{equation}
where $|\Omega|_i$ is the pullback of $|\Omega|:=|\Omega/\!\vol\!|$ to
$\gamma^1$ with
its new metric. Note also that pulling functions up from $\gamma$ and
taking their exterior derivative $d$ on either $\gamma^1$ or on
$\gamma^n$ gives the same result via the obvious inclusion
$\gamma^1\subset\gamma^n$ commuting with the projections to $\gamma$.
Again we want to control this evolution equation as $i\to\infty$.

$|d\theta|=|\theta_s|=|\partial_s[\theta(\gamma^1_s)/2+(n/2-1)
\theta(p(\gamma^1))]|$, so this is bounded by the estimates above,
for $i$ sufficiently large that $\theta(p(\gamma^1))$ is $C^1$
close to $\theta^1/2+$\,const. in the disc in which we are working.
So we can bound the last term in (\ref{eqn}) by a constant times
$$
\kappa_i^{-1}{\sup|d\Omega|\over\inf|\Omega|},
$$
where the sup and inf are taken in the original metric over a small
neighbourhood of $(0,\ldots,0,z)\in X^n$. This tends to zero as
$i\to\infty$.

Computing the Laplacian on the space $\gamma^n$, with a radial
coordinate $s$ and rotational symmetry about the origin $s=0$,
makes (\ref{eqn})
\begin{equation} \label{tat}
\theta_t=\theta_{ss}+(n-1){R^i_s\over R^i}\theta_s+O(\kappa_i^{-1}),
\end{equation}
where $R^i=R^i(s)$ is the radius of the sphere $S^{n-1}$ at $s$ in the
new metric.

To compute $R^i$, we use our $\gamma^1\subset X^1$ arclength coordinate
$s$, the radial coordinate $r$ on $X^1$, and a radial coordinate
$\rho$ on $\C$. For $i$ sufficiently large, for $s\le\kappa_i^{1/2}$
in the new metric, we can approximate $p$ linearly about $z$ and
so assume that $R^i$ is as close as we like
to $\kappa_i\sqrt{|\dot p(z)|\rho}$ in $C^2$. Therefore $R^i_r$ is
approximated by
\begin{equation} \label{grig}
R^i_r=\frac{R^i_\rho}{r_\rho}=\frac{R^i_\rho}{\sqrt{\kappa_i^2+(R^i_\rho)^2}}
=\frac{1}{\sqrt{1+\frac{\kappa_i^2}{(R^i_\rho)^2}}}\approx
\frac{1}{\sqrt{1+\frac{4\rho}{|\dot{p}(z)|}}}\ ,
\end{equation}
which is bounded and tends to $1$ in the interval
$s\in[0,\kappa_i^{1/2})$. Similarly $R^i_{rr}$ can be taken to be
arbitrarily small on the same interval, for $i$ sufficiently large.

Since $r_s$ is bounded, this gives bounds on $R^i_s$, implying that,
on passing to a subsequence if necessary, the functions $R^i(s)$
are convergent as $i\to\infty$ by the Arzel\`a-Ascoli theorem.

Note also that for all $i$, $R^i_s(0)=1$. But to preserve this
in the limit, we must similarly bound $R^i_{ss}$. Differentiating
$r^2_s+r^2(\theta^1_s)^2=1$ and $R^i_s=R^i_rr_s$ gives
\begin{equation} \label{final}
R^i_{ss}=R^i_{rr}r_s^2-R^i_r(r(\theta^1_s)^2+
r^2\theta^1_s\theta^1_{ss}/r_s).
\end{equation}
We have bounded $R^i_{rr},\ r_s,\ r_s^{-1},\ r\theta^1_s$ and
$\theta^1_s$; all this leaves is the last term in (\ref{final}).

We have approximated $\theta^1_s$ in $C^1$ by
${1\over r}\sin(\theta(\gamma^1_s)-\theta^1)$; differentiating
approximates $r^2\theta^1_s\theta^1_{ss}/r_s$ as closely as
we like (as $i\to\infty$) to
$$
{r\theta^1_s\over r_s}(\kappa_i-\theta^1_s)\cos(\theta(\gamma^1_s)
-\theta^1)-\theta^1_s\sin(\theta(\gamma^1_s)-\theta^1),
$$
which we have bounded already.
In conclusion, after passing to a subsequence if necessary, $R^i$ is
$C^1$ convergent to some $R$ with $R_s(0)=1$, and the phase function
$\theta^\infty$ of the limit curve $\gamma^1_\infty$
(which exists by Arzel\`a-Ascoli
since $|\gamma^1_s|=1$ and $|\gamma^1_t|$ is bounded by the bound on
its curvature $\kappa$) satisfies the limit of (\ref{tat}):
\begin{equation} \label{lap}
\theta^\infty_t=\theta^\infty_{ss}+(n-1){R_s\over R}\theta^\infty_s.
\end{equation}

But this is just the heat equation for $\theta$ on $\R^n$ with radial
coordinate $s$ and the $O(n)$-invariant metric in which the radius of
the $S^{n-1}$ fibre over $s$ is $R(s)$. By construction of the time
rescaling (\ref{scale}) it exists for all time $t\le0$, and the
solution $\theta^\infty$ is bounded. Also the metric is uniformly
elliptic compared to the flat metric, as $R_r=1$ (\ref{grig}) and
$r_s$ is bounded above and below away from zero. Therefore, by
Moser's Harnack inequality \cite{Mo}, $\theta^\infty$ is in fact
constant.

But for all $i$, max\,$\theta_s=1$ by construction, and passing to
a subsequence if necessary the point
where the maximum is obtained is convergent. To show then that
max\,$\theta^\infty_s=1$, to get our contradiction, we
need only know that $\theta_s$ is, say, \emph{uniformly (in $i$)}
H\"older continuous. This is again a consequence of (\cite{LSU}
Section III Theorem 10.1) as follows.
By the boundedness of $\theta_s$ and the parabolic
equation it satisfies (different for each $i$), $\theta_s$ is
in fact $\alpha$-H\"older continuous for some $\alpha>0$, with
$H^\alpha$ norm bounded by the bounds on the coefficients of
the parabolic equations. But these are bounded uniformly in $i$,
as a glance at (\ref{eqn}) confirms: the correction term tends to
zero, and the Laplacian term (and its derivative with respect to
$s$) is controlled by $C^2$ bounds on the metric which we provided
above by bounding $R^i(s),\,R^i_s$ and $R^i_{ss}$.
\end{Proof}

Finally we show that this infinite time flow converges using
standard techniques (see for instance \cite{C} for a harder result).
Notice that the same scaling proof (\ref{blowup}) that the
curvature of $\gamma^1$ does not blow up in finite time shows the
same for our now infinite time flow. So, using the $O(n)$ symmetry,
the curvature of the metric on $\gamma^n$ stays uniformly bounded,
and we have a $C^1$ bound on $\theta$. Therefore in the equation
(\ref{newdot}) for $\theta$ on $\gamma^n$, which we rewrite as
\begin{equation} \label{ntd}
\dot\theta=-\Delta^\Omega\theta:=|\Omega|^{-1}*d(|\Omega|*d\theta),
\end{equation}
the coefficients have at least uniform $C^1$ bounds; $\theta$ then
acquires a uniform $C^3$ bound by parabolic Schauder estimates (see
\cite{LSU} III Theorem 12.1, for instance). Again by $O(n)$ symmetry
we now get uniform $C^2$ bounds on $\gamma^n$'s curvature. And so it
goes on, giving $C^\infty$ bounds and allowing us to extract a
subsequence of times for which the flow converges in $C^\infty$ to
a Slag.

To see that the flow converges without having to pass to
a subsequence we need only show convergence of $\theta$ in $L^2$;
this way no other subsequence of the flow can converge to
a different limit in $C^\infty$. In fact we use an $L^2$-norm
weighted by $|\Omega|$, and compute using (\ref{ntd})
and (\ref{nvdot}):
$$
{d\over dt}\int_{\gamma^n}(\theta-\bar\theta)^2|\Omega|\vol=
\int_{\gamma^n}\left\{2(\theta-\bar\theta)(-\Delta^\Omega\theta-
{d\over dt}\bar\theta)-(\theta-\bar\theta)^2|d\theta|^2\right\}
|\Omega|\vol,
$$
where $\bar\theta=\int\theta|\Omega|\vol/\int|\Omega|\vol$ is
constant on $\gamma^n$ (but not in time). This is then bounded
above by
$$
-2\int(\theta-\bar\theta)\Delta^\Omega(\theta-\bar\theta)
|\Omega|\vol.
$$
It is easily checked that $\Delta^\Omega=d_{\,}^{*^\Omega}
d$, where $d_{\,}^{*^\Omega}$ is the adjoint of $d$ with
respect to the $L^2$-metric $\int_{\gamma^n}\langle\ \cdot\ ,
\ \cdot\ \rangle|\Omega|\vol$ we are using. So its kernel is just
the constants, and by the uniform $C^\infty$ bounds on the metric
of $\gamma^n$ and $|\Omega|$ we can get a uniform lower bound
$\lambda>0$ for its first nonzero eigenvalue. This then gives
a bound $\int(\mu\Delta\mu)|\Omega|\vol\ge\lambda\int\mu^2|\Omega|\vol$
for functions $\mu$ of integral zero. Setting $\mu=\theta-\bar\theta$
gives
$$
{d\over dt}\int(\theta-\bar\theta)^2|\Omega|\vol\le-2\lambda
\int(\theta-\bar\theta)^2|\Omega|\vol,
$$
which therefore tends to zero as required. Theorem \ref{faq} is proved. \\

We end by noting that one can get cone-type bounds similar to
those of Lemma \ref{cone} on the tangent
direction $\theta(\gamma^1_s)$ even as a curve $\gamma$
representing an \emph{unstable} Lagrangian $\gamma^n$ approaches
and breaks across a zero of $p$. Suppose that the initial phase
variation of $\gamma^n$ is, without loss of generality,
in some $(-\delta,\delta)$. Draw the cone with boundary the
SLags emanating from the zero of $p$ with phase $-\delta,\ \pi+
\delta$ (the straight lines emanating from the zeros of $p$ drawn
in Figures 2 and 4 display the $\delta=0$ cone for dimensions 2
and 3 respectively). Then in a
sufficiently small neighbourhood of the zero, the variation in
$\theta(p)$ can be taken to be less than $2\pi/n+2\delta+\epsilon$
for any $\epsilon>0$. Since we are free to make $\delta$ slightly
smaller without violating the initial bounds on phase, we can
ensure that there is a neighbourhood of the zero of $p$, and a cone
with vertex at $p$ whose walls $\gamma$ cannot cross, such that for
the part of $\gamma$ lying in this neighbourhood, the variation
$$
\sup\theta(p(\gamma))-\inf\theta(p(\gamma))<2\pi/n+2\delta. 
$$
Comparing with (\ref{pf}) gives
$$
\sup\theta(\gamma')-\inf\theta(\gamma')<(n/2-1)(2\pi/n+2\delta)+
2\delta=(1-2/n)\pi+n\delta
$$
so that again no $180^o$ kinks can occur while this is less than
$\pi$, i.e. for $\delta\le2\pi/n^2$. So again the analysis should
be tractable in this case (more general spiralling around a zero
would make matters worse). However, we have not carried out the
analysis necessary to show that at the moment $\gamma$ reaches
the zero of $p$ it is sufficiently smooth that the two resulting
curves it splits into give $C^2$ Lagrangians (whose flow we could
restart). \\

Of course by just studying the simple examples above we cannot
hope to know how bad the singularities are that arise in finite
time in the general case.
Also, as mentioned in \cite{Th}, we should perhaps restrict to
those Lagrangians whose Floer cohomology is well defined \cite{FO3}.
This includes all homology spheres, however.

We should also point out the obvious fact that most of the
evidence for our conjecture, other than perhaps the mirror
symmetry and study of Joyce's examples in \cite{Th}, has been
essentially one-dimensional (either
for $T^2$ in \cite{Th}, or by symmetry reduction in the examples
above). This is unrepresentative, essentially because the angles at
which Lagrangians intersect (the $\alpha_i$s of (\ref{l1}))
are all the same in this situation, and so are determined by the
phase (their sum). So interesting phenomena, where degrees in
Floer cohomology change (e.g. a Hom becomes an Ext$^i$ on the mirror
while the phase remains fixed) are largely lost due to them being
controlled entirely by the phase.

\small \noindent {\tt richard.thomas@ic.ac.uk} \newline
\noindent Department of Mathematics, Imperial College, Huxley
Building, 180 Queen's Gate, London, SW7 2BZ. UK. \newline\newline
\noindent {\tt yau@math.harvard.edu} \newline
\noindent Department of Mathematics, Harvard University, One
Oxford Street, Cambridge MA 02138. USA.

\end{document}